\documentclass[12pt]{article}

\usepackage{graphics,graphicx}
\usepackage{subfigure}

\usepackage{amssymb, amsfonts, amsbsy, amsmath,mathrsfs} 
\usepackage{eucal}

\usepackage{multirow}
\usepackage{lscape}
\usepackage{setspace}

\newcommand{\Eref}[1]{Equation (\ref{#1})}
\newcommand{\Erefs}[1]{Equations (\ref{#1})}

\newcommand{\fref}[1]{Figure~\ref{#1}}
\newcommand{\frefs}[1]{Figures~\ref{#1}}

\newcommand{\ff}{\mathbf{f}}
\newcommand{\KK}{\mathbf{K}}
\newcommand{\bm}{\mathbf{M}}

\newcommand{\bveps}{\boldsymbol{\varepsilon}}
\newcommand{\bvsig}{\boldsymbol{\sigma}}

\begin{document}

\begin{center}
\large Bending and vibration of functionally graded material sandwich plates using an accurate theory\footnote{Preprint sent to Finite Elements in Analysis and Design}
\end{center}

\begin{center}{S Natarajan$^{a,\dagger}$, M Ganapathi$^{b}$} 
\end{center}

\begin{center}\small{
$^{a,\dagger}$Institute of Mechanics and Advanced Materials, Theoretical and Computational Mechanics, Cardiff University, UK. \\ Email: sundararajan.natarajan@gmail.com  \\
$^{b}$Head, Stress \& DTA, IES-Aerospace, Mahindra Satyam Computers Services Ltd., Bangalore, India\\}
\end{center}

\begin{abstract}
In this paper, the bending and the free flexural vibration behaviour of sandwich functionally graded material (FGM) plates are investigated using QUAD-8 shear flexible element developed based on higher order structural theory. This theory accounts for the realistic variation of the displacements through the thickness. The governing equations obtained here are solved for static analysis considering two types of sandwich FGM plates, viz., homogeneous face sheets with FGM core and FGM face sheets with homogeneous hard core. The in-plane and rotary inertia terms are considered for vibration studies. The accuracy of the present formulation is tested considering the problems for which three-dimensional elasticity solutions are available. A detailed numerical study is carried out based on various higher-order models to examine the influence of the gradient index and the plate aspect ratio on the global/local response of different sandwich FGM plates.

\end{abstract}

\begin{footnotesize}
\textbf{Keywords}: Functionally graded material plate, higher-order theory, thermal loading, mechanical loading, free vibration, shear flexible element.
\end{footnotesize}


\section{Introduction}
With ever-increasing demand for high strength-to-weight ratio materials, the engineered materials have been replacing the conventional materials in automotive, nuclear and aerospace industries. In general, most of the engineered materials are inspired from nature. A new class of materials was introduced by Japan scientists~\cite{Koizumi1993} to decrease the thermal stresses in the propulsion systems and the airframe of the space planes. This class of engineered materials was coined as functionally graded material (FGM). These materials are made up of mixture of ceramics and metals, that are characterized by the \emph{smooth and continuous variation} in the properties from one surface to another~\cite{Koizumi1993,Reddy2000}. For the structural integrity, the FGMs are preferred over the fiber-matrix composites that may results in debonding due to the mismatch in the mechanical properties across the interface of two discrete materials bonded together. This has attracted lot of researchers in understanding the mechanics and mechanism of the FGM structures.

Among the various structural constructions, the sandwich types of structures are commonly used in the aerospace vehicles, because of its outstanding bending rigidity, low specific weight, excellent vibration characteristics and good fatigue properties. Laminated composite types of constructions are, in general, adopted in sandwich structures. However, due to the sudden change in the material properties from one layer to another, the variation of the interfacial stress distribution is significant at the facesheet-core interface. Furthermore, the response of such laminated composites depend on the bonding characteristics. In contrast, the FGM sandwich can alleviate the large variation in the interfacial stress distribution, because of the gradual variation of the material properties at the facesheet-core interface. For predicting the realistic structural behaviour of such components, more accurate analytical/numerical analysis based on the three-dimensional models may be computationally involved and expensive. Hence, among the researchers, there is a growing appreciation of the importance of applying two-dimensional theories with new kinematics for the evolution of the accurate structural analysis. Various structural theories proposed for the FGM structures have been examined and some of the important contributions pertaining to the sandwich FGM plates are discussed here. The most commonly used FGM sandwich constructions are: the  FGM facesheet with homogeneous core and the homogeneous facesheet with  FGM core. These sandwich constructions can be considered for the requirement of light weight and high bending stiffness in design (for example, the placement of actuators/sensors in the field control and thermal/mechanical load bearing component design) by appropriately selecting the soft/hard core metal or the ceramic layers.

Venkataraman and Sankar~\cite{Venkataraman2001} and Anderson~\cite{Anderson2003} have studied the effect of FGM core in sandwich beam on the interfacial shear stresses.  Pan and Han~\cite{Pan2005} analysed the static response of the multilayered rectangular plate made of functionally graded, anisotropic and linear magneto-electro-elastic materials. Das \textit{et al.,}~\cite{Das2006} have investigated a sandwich consisting of FGM soft core with relatively stiff orthotropic facesheets employing triangular plate element. Ganesan~\textit{et al.,}~\cite{Bhangale2006} studied static and dynamic response of FGM plate with viscoelastic core, whereas Shen~\cite{Shen2005} examined FGM sandwich plates with piezo-electric core subjected to thermo-electro-mechanical loading. Zenkour~\cite{Zenkour2005,Zenkour2005a} studied analytically the static and dynamics of sandwich FGM plates with homogeneous hard core using sinusoidal shear deformation theory.  Li et al. ~\cite{Li2008} have presented three-dimensional analytical solutions for multi-layer FGM sandwich plates based on Ritz method in conjunction with Chebhyshev polynomial series. It is observed from these studies that first- order, third-order, and sinusoidal shear deformation theories have been extensively used for the analysis of sandwich FGM plates.  However to author's knowledge, the theories accounting the variation of in-plane displacement through the thickness,  and the possible discontinuity in slope at the interface, and the thickness stretch affecting transverse deflection is not exploited while investigating the structural behaviour of FGM sandwich structures. A Layer wise theory is the possible candidature for this purpose, but it may be computationally expensive as the number of unknowns to be solved increases with the increase in the number of mathematical or physical layers. Ali et al.~\cite{Ali1999}, and Ganapathi and Makhecha~\cite{Ganapathi2001} have employed a new higher-order plate theory based on global approach, for multi-layered laminated composites by incorporating the realistic through the thickness approximations of the in-plane and transverse displacements by adding a zig-zag functions and higher-order terms, respectively. This formulation has proved to give very accurate results for the composite laminates. Such model for the current problems is now explored as a candidature while comparing with the three-dimensional model.

In this paper, a $\mathcal{C}^o$ 8-noded quadrilateral plate element with 13 degrees of freedom per node ~\cite{Ganapathi2001,Makhecha2001} based on the higher order theory~\cite{Ali1999} is employed to study the static deflection and the free vibration analysis of thick/thin sandwich functionally graded material plates. The efficacy of the present formulation, for the static analyses subjected to the thermal/mechanical loads and the free vibration analysis is illustrated through the numerical studies. The accuracy of the present element with those of the other two-dimensional numerical/analytical models and the three-dimensional elasticity analysis are discussed considering various parameters.

The paper is organized as follows. The next section will give an introduction to the FGM and a brief overview of the higher order theory. Section \S \ref{eledes} describes the 8-noded quadrilateral plate element used in the current study. Section \S \ref{numexamples} present numerical results for the static deflection and the free vibration of thick/thin sandwich FGM plates, followed by concluding remarks in Section \S \ref{conclu}.

\section{Theoretical Formulation}
\label{formulate}

\subsection{Functionally graded material plate}
\begin{figure}[htpb]
\centering
\input{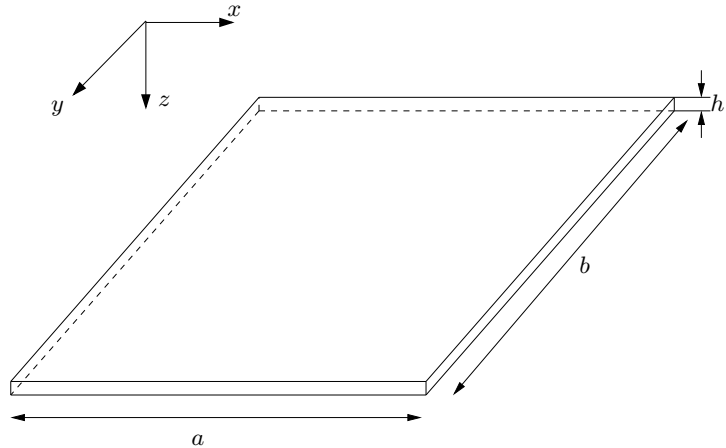}
\caption{Co-ordinate system of rectangular FGM plate. $x,y$ along the in-plane directions and $z$ along the thickness direction.}
\label{fig:fgmplate}
\end{figure}

Consider a rectangular FGM plate with co-ordinates $x,y$ along the in-plane directions and $z$ along the thickness direction as shown in the \fref{fig:fgmplate}. The material is assumed to be graded only in the thickness direction according to a power-law distribution. The homogenized material properties are computed using the rule of mixtures.

\subsubsection*{Estimation of mechanical properties}
The effective material properties for each layer, viz., Young's modulus, Poisson's ratio and the mass density is estimated by the following power-law function:

\begin{equation}
P(z) = P_m V_m + P_c V_c
\label{eqn:powerlaw}
\end{equation}

Here, $V_i~(i=c,m)$ is the volume fraction of the phase material. The subscripts $c$ and $m$ refer to the ceramic and the metal phases, respectively. Note that $P_m$ and $P_c$ are the properties of the metallic and the ceramic phases, respectively. The volume fractions of the ceramic and the metal phases are related by $V_c + V_m = 1$, and $V_c$ is expressed as:

\begin{equation}
V_c(z) = \left( {2z + h \over 2h} \right)^n, \hspace{0.2cm}  n \ge 0
\label{eqn:volFrac}
\end{equation}

where $n$ in \Eref{eqn:volFrac} is the volume fraction exponent, also referred to as the gradient index in the literature. The properties of the FGM plate vary continuously through the thickness based on a power-law function, given by~\Eref{eqn:powerlaw}. In this study, the following two types of power-law FGMs are considered.

\subsubsection*{Type A - FGM face sheet and homogeneous hard core}
In this case, the volume fraction of the FGMs is assumed to follow 

\begin{eqnarray}
\renewcommand{\arraystretch}{2}
V_1(z) = \left( \frac{z - z_1}{z_2 - z_1} \right)^n, \hspace{1cm} z \in [z_1,z_2] \nonumber \\
V_2(z) = 1, \hspace{1cm} z \in [z_2,z_3] \nonumber \\
V_3(z) = \left( \frac{z - z_4}{z_3 - z_4} \right)^n, \hspace{1cm} z \in [z_3,z_4]
\label{eqn:TypeA_volFrac}
\end{eqnarray}

where $V_i (i=1,2,3)$ is the volume fraction of layer $i$ and $n$ is the gradient index. The core is considered as a fully ceramic material. The top and bottom surfaces of the plate are metal-rich.

\subsubsection*{Type B - Homogeneous face sheet and FGM core}
In this case, the volume fraction of the FGMs is assumed to follow 

\begin{eqnarray}
\renewcommand{\arraystretch}{2}
V_1(z) = 0, \hspace{2cm} z \in [z_1,z_2] \nonumber \\
V_2(z) = \left( \frac{z - z_2}{z_3 - z_2} \right)^n, z \in [z_2,z_3] \nonumber \\
V_3(z) = 1, \hspace{2cm} z \in [z_3,z_4]
\label{eqn:TypeB_volFrac}
\end{eqnarray}

\begin{figure}
\subfigure[FGM 1-1-1]{\includegraphics[scale=0.45]{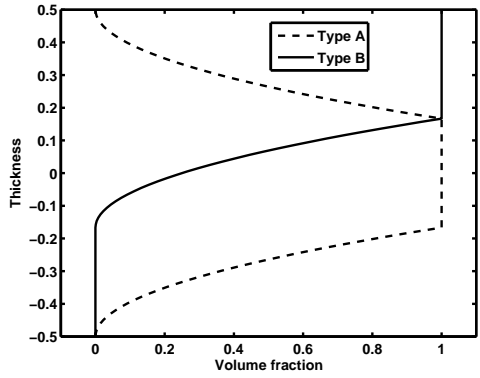}}
\subfigure[FGM 1-2-1]{\includegraphics[scale=0.45]{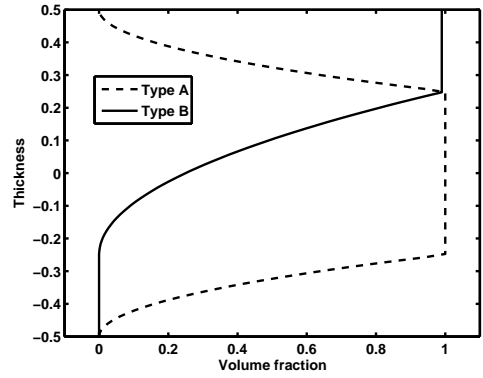}}
\subfigure[FGM 2-2-1]{\includegraphics[scale=0.45]{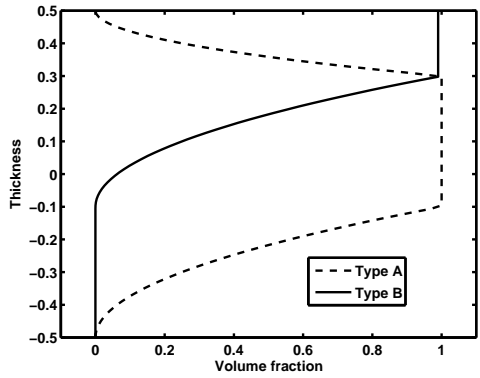}}
\subfigure[FGM 1-8-1]{\includegraphics[scale=0.45]{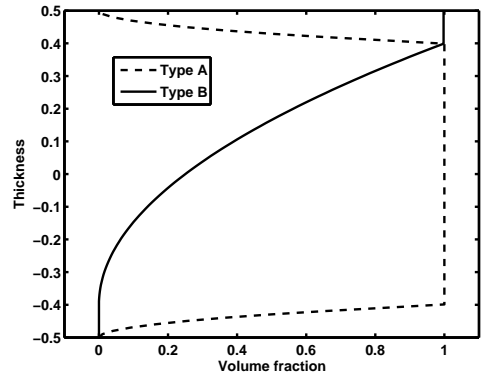}}
\caption{Through thickness variation of volume fraction for various types of sandwich FGM plates for gradient index $n=$2.}
\label{fig:volfrac}
\end{figure}

The core is considered as a functionally graded material. The top surface of the plate is metal rich and the bottom surface of the plate is ceramic rich. \fref{fig:volfrac} shows the variation of the volume fractions of ceramic and metal, respectively, in the thickness direction $z$ for the two types of sandwich FGM plates considered in this study.

\subsection{Plate formulation}
The sandwich FGM plate is assumed to be made up of three discrete layers. The material properties for $k^{th}$ layer is governed by \Eref{eqn:TypeA_volFrac} or \Eref{eqn:TypeB_volFrac}. The in-plane displacements $u^k$ and $v^k$, and the transverse displacement $w^k$ for the $k^{th}$ layer, are assumed as~\cite{Ali1999,Makhecha2001,Ganapathi2001}:

\begin{eqnarray}
\begin{split}
u^k(x,y,z,t) = u_o(x,y,z,t) + z \theta_x(x,y,t) + z^2 \beta_x(x,y,t) + \\ z^3 \phi_x(x,y,t) + S^k \psi_x(x,y,t) \end{split} \nonumber \\
\begin{split}
v^k(x,y,z,t) = v_o(x,y,z,t) + z \theta_y(x,y,t) + z^2 \beta_y(x,y,t) + \\ z^3 \phi_y(x,y,t) + S^k \psi_y(x,y,t) \end{split} \nonumber \\
w^k(x,y,z,t) = w_o(x,y,t) + z w_1(x,y,t) + z^2 \Gamma(x,y,t)
\label{eqn:dispField}
\end{eqnarray}

The terms with even powers of $z$ in the in-plane displacements and odd powers of $z$ occurring in the expansion for $w^k$ correspond to the stretching problem. However, the terms with odd powers of $z$ in the in-plane displacements and the even ones in the expression for $w^k$ represent the flexure problem. $u_o, v_o$ and $w_o$ are the displacements of a generic point on the reference surface; $\theta_x$ and $\theta_y$ are the rotations of the normal to the reference surface about the $y$ and $x$ axes, respectively; $w_1,\beta_x,\beta_y,\Gamma,\phi_x$ and $\phi_y$ are the higher-order terms in the Taylor`s series expansions, defined at the reference surface. $\psi_x$ and $\psi_y$ are generalized variables associated with the zigzag function, $S^k$. The zigzag function, $S^k$, as given in~\cite{Murukami1986,Ganapathi2001,Roderigues2011}, is defined by

\begin{equation}
S^k = 2(-1)^k \frac{z_k}{h_k}
\end{equation}

where $z_k$ is the local transverse coordinate with the origin at the center of the $k^{th}$ layer and $h_k$ is the corresponding layer thickness. Thus, the zigzag function is piecewise linear with values of -1 and 1 alternatively at different interfaces. The `zig-zag' function, as defined above, takes care of the inclusion of the slope discontinuities of $u$ and $v$ at the interfaces of the sandwich plate as observed in the exact three-dimensional elasticity solutions of thick sandwich functionally graded materials. The main advantage of using such formulation is that such a function is more economical than a discrete layer approach~\cite{Nosier1993,Ferreira2005}.

The strains in terms of mid-plane deformation, rotations of normal and higher order terms associated with displacements are:

\begin{equation}
\bveps = \left \{ \begin{array}{c} \bveps_{\rm bm} \\ \bveps_s \end{array} \right\}
\label{eqn:strain}
\end{equation}

The vector  $\bveps_{\rm bm}$ includes the bending and membrane terms of the strain components and vector $\bveps_s$ contains the transverse shear strain terms. These strain vectors are defined as ~\cite{Ali1999,Ganapathi2001,Makhecha2001}:

\begin{eqnarray}
\bveps_{\rm bm} &=& \left\{ \begin{array}{c} \varepsilon_{xx} \\ \varepsilon_{yy} \\ \varepsilon_{zz} \\ \gamma_{xy} \end{array} \right\} + \left\{ \begin{array}{c} u_{,x} \\ v_{,y} \\ w_{,z} \\ u_{,y} + v_{,x} \end{array} \right\} \nonumber \\
&=& \bveps_0 + z \bveps_1 + z^2 \bveps_2 + z^3 \bveps_3 + S^k \bveps_4 
\end{eqnarray}

\begin{eqnarray}
\bveps_s &=& \left\{ \begin{array}{c} \gamma_{xz} \\ \gamma_{yz} \end{array} \right\} = \left\{ \begin{array}{c} u_{,z} + w_{,x} \\ v_{,z} + w_{,y} \end{array} \right\} \nonumber \\
&=& \gamma_o + z \gamma_1 + z^2 \gamma_2 + S^k_{,z} \gamma_3
\end{eqnarray}

where,

\begin{eqnarray}
\bveps_o = \left\{ \begin{array}{c} u_{o,x} \\ v_{o,y} \\ w_1 \\ u_{o,y} + v_{o,x} \end{array} \right\},\hspace{1cm} \bveps_1 = \left\{ \begin{array}{c} \theta_{x,x} \\ \theta_{y,y} \\ 2\Gamma \\ \theta_{x,y} + \theta_{y,x} \end{array} \right\}, \nonumber \\
\bveps_2 = \left\{ \begin{array}{c} \beta_{x,x} \\ \beta_{y,y} \\ 0 \\ \beta_{x,y} + \beta_{y,x} \end{array} \right\}, \hspace{1cm} \bveps_3 = \left\{ \begin{array}{c} \phi_{x,x} \\ \phi_{y,y} \\ 0 \\ \phi_{x,y} + \phi_{y,x} \end{array} \right\}, \nonumber \\
\bveps_4 = \left\{ \begin{array}{c} \psi_{x,x} \\ \psi_{y,y} \\ 0 \\ \psi_{x,y} + \psi_{y,x} \end{array} \right\}.
\end{eqnarray}

and,

\begin{eqnarray}
\gamma_o = \left\{ \begin{array}{c} \theta_x + w_{o,x} \\ \theta_y + w_{o,y} \end{array} \right\}, \hspace{1cm} \gamma_1 = \left\{ \begin{array}{c} 2\beta_x + w_{1,x} \\ 2\beta_y + w_{1,y} \end{array} \right\}, \nonumber \\
\gamma_2 = \left\{ \begin{array}{c} 3\phi_x + \Gamma_{,x} \\ 3\phi_y + \Gamma_{,y} \end{array} \right\}, \hspace{1cm} \gamma_3 = \left\{ \begin{array}{c} \psi_x S_{,z}^k \\ \psi_y S_{,z}^k \end{array} \right\} .
\end{eqnarray}

The subscript comma denotes partial derivatives with respect to the spatial coordinate succeeding it. The constitutive relations for an arbitrary layer $k$ can be expressed as:

\begin{eqnarray}
\bvsig &=& \left\{ \begin{array}{ccccccc} \sigma_{xx} & \sigma_{yy} & \sigma_{zz} & \sigma_{xy} & \sigma_{xz} & \sigma_{yz} \end{array} \right\}^{\rm T} \nonumber \\
&=& \bf{Q}^k \left\{ \begin{array}{cc} \bveps_{\rm bm} & \bveps_s \end{array} \right\}^{\rm T}
\end{eqnarray}

where $\bf{Q}_k$ is the stiffness coefficient matrix defined as:

\begin{eqnarray}
Q_{11}^k = Q_{22}^k = {E(z) \over 1-\nu^2}; \hspace{1cm} Q_{12}^k = {\nu(z) E(z) \over 1-\nu(z)^2}; \hspace{1cm} Q_{16}^k = Q_{26}^k = 0 \nonumber \\
Q_{44}^k = Q_{55}^k = Q_{66}^k = {E(z) \over 2(1+\nu(z)) }
\label{eqn:stiffcoeff}
\end{eqnarray}

where the modulus of elasticity $E(z)$ and Poisson's ratio $\nu(z)$ are given by \Eref{eqn:powerlaw}. The governing equations of motion are obtained by applying Lagrangian equations of motion given by:

\begin{equation}
\frac{d}{dt} \left[ \frac{\partial (T-U)}{\partial \dot{\delta_i}} \right] - \left[ \frac{\partial (T-U)}{\partial \delta_i} \right] = 0, \hspace{1cm} i = 1,2,\cdots,n
\label{eqn:lagrange}
\end{equation}

where $\delta_i$ is the vector of degrees of freedom and $T$ is the kinetic energy of the plate given by:

\begin{equation}
T(\delta) = \frac{1}{2} \iint \left[ \sum_{k=1}^n \int\limits_{h_k}^{h_{k+1}} \rho_k \left\{ \begin{array}{ccc} \dot{u}_k & \dot{v}_k & \dot{w}_k \end{array} \right\} \left\{ \begin{array}{ccc} \dot{u}_k & \dot{v}_k & \dot{w}_k \end{array}\right\} ^{\rm T}~dz \right] dx dy
\label{eqn:kinetic}
\end{equation}

where $\rho_k$ is the mass density of the $k^{th}$ layer, $h_k$ and $h_{k+1}$ are the $z$ coordinates corresponding to the bottom and top surfaces of the $k^{th}$ layer. The total potential energy function $U$ is given by:

\begin{equation}
U(\delta) = \frac{1}{2} \iint \left[ \sum_{k=1}^n \int\limits_{h_k}^{h_{k+1}} \bvsig^{\rm T} \bveps ~dz \right] dx dy - \iint {\bf q} w~dxdy
\label{eqn:potential}
\end{equation}

where ${\bf q}$ is the distributed force acting on the top surface of the plate. Substituting \Erefs{eqn:potential} and (\ref{eqn:kinetic}) in \Eref{eqn:lagrange}, one obtains the following governing equations for static deflection and free vibration of plate.

\paragraph*{Static deflection}
\begin{equation}
\KK \boldsymbol{\delta} = \ff
\label{eqn:staticdefl}
\end{equation}

\paragraph*{Free vibration}
\begin{equation}
\bm \ddot{\boldsymbol{\delta}} + \KK \boldsymbol{\delta} = \bf{0}
\label{eqn:freevib}
\end{equation}

where $\bm$ is the mass matrix, $\KK$ is the stiffness matrix and $\ff$ is the external force vector. In the present study, while performing the integration, terms having thickness co-ordinate $z$ are integrated with higher order Gaussian quadrature, because the material properties vary continuously through the thickness. The terms containing $x$ and $y$ are evaluated using full integration with 3 $\times$3 Gauss integration rule. The frequencies and mode shapes are obtained from \Eref{eqn:freevib} using the standard generalized eigenvalue algorithm.

\section{Element description}
\label{eledes}
In this paper, $\mathcal{C}^o$ continuous, eight-noded serendipity quadrilateral shear flexible plate element with 13 nodal degrees of freedom \\$(u_o,v_o,w_o,\theta_x,\theta_y,w_1,\beta_x,\beta_y,\Gamma,\phi_x,\phi_y,\psi_x,\psi_y:$ 13 DOF is used$)$. The finite element represented as per the kinematics based on \Eref{eqn:dispField} is referred to as Q8-HSDT13 with cubic variation. Interested readers are referred to the literature~\cite{Ganapathi2001,Makhecha2001}, where the behaviour of the element is described in detail. The element is shown to be free from locking syndrome, absence of spurious energy modes, passes patch test and exhibits faster convergence~\cite{Ganapathi2001,Makhecha2001}.
Three more alternate discrete models are proposed to study the influence of higher-order terms in the displacement functions, whose displacement fields are deduced from the original element by deleting the appropriate degrees of freedom. These alternate models, and the corresponding degrees of freedom are listed in Table \ref{table:alternatemodels}

\begin{table} [htpb]
\renewcommand\arraystretch{1.5}
\caption{Alternate eight-noded finite element models}
\centering
\begin{tabular}{ll}
\hline
Finite element model & Degrees of freedom per node  \\
\hline
HSDT13 (present) & $u_o,v_o,w_o,\theta_x,\theta_y,w_1,\beta_x,\beta_y,\Gamma,\phi_x,\phi_y,\psi_x,\psi_y$ \\
HSDT11 & $u_o,v_o,w_o,\theta_x,\theta_y,w_1,\beta_x,\beta_y,\Gamma,\phi_x,\phi_y$ \\
HSDT9 & $u_o,v_o,w_o,\theta_x,\theta_y,\beta_x,\beta_y,\phi_x,\phi_y$ \\
FSDT & $u_o,v_o,w_o,\theta_x,\theta_y$ \\
\hline
\end{tabular}
\label{table:alternatemodels}
\end{table}

\section{Numerical results and discussion}
\label{numexamples}
In this section, we present the static response and the natural frequencies of sandwich FGM plates using the eight-noded shear flexible quadrilateral element. The effect of plate aspect ratio and gradient index are studied. In this study, only simply supported boundary conditions are considered and are as follows:

\begin{eqnarray}
u_o = w_o = \theta_x = w_1 = \Gamma = \beta_x = \phi_x = \psi_x = 0, \hspace{0.2cm} ~\textup{on} ~ y = 0,b \nonumber \\
v_o = w_o = \theta_y = w_1 = \Gamma = \beta_y = \phi_y = \psi_y = 0, \hspace{0.2cm} ~\textup{on} ~ x= 0,a
\end{eqnarray}

where $a$ and $b$ refer to the length and width of the plate, respectively. The FGM plate considered here consists of Alumina and Aluminum. The mass density $\rho$ and Young's modulus $E$ are: $\rho_c = $ 3800 kg/m$^3$, $E_c=$ 380 GPa for Alumina and $\rho_m=$ 2707 kg/m$^3$, $E_m=$ 70 GPa for Aluminum. Poisson's ratio $\nu$ is assumed to be constant and taken as 0.3 for the current study. Here, the modified shear correction factor obtained based on energy equivalence principle~\cite{Singh2011,Natarajan2011} is used for the FSDT model. The transverse shear stresses are evaluated by integrating the three-dimensional equilibrium equations for all types of elements. For the current study, three different core thickness (1-1-1, 1-2-1, 2-2-1) for both Type A and Type B sandwich FGM plates, three thickness ratio $a/h$ (5,10,100) and four gradient indices $n$ (0,0.5,1,5) are considered.

\subsection{Static analysis}

The static analysis is conducted for Type A FGM sandwich plate. The following two types of loading are considered:

\begin{itemize}
\item Mechanical loading: $q = q_o \sin (\pi x/a) \sin (\pi y/b)$,
\item Thermal loading: $T = T_o (2z/h) \sin(\pi x/a) \sin( \pi y/b)$.
\end{itemize}

where $q_o$ and $T_o$ are the amplitudes of mechanical and thermal loads, respectively. The physical quantities are nondimensionalized by relations:

\begin{eqnarray}
\renewcommand{\arraystretch}{2}
(\overline{u},\overline{v}) = \frac{100E_o}{q_o h S^3}(u,v) \nonumber \\
\overline{w} = \frac{100E_o}{q_o h S^4}(w) \nonumber \\
(\overline{\sigma}_{xx},\overline{\sigma}_{yy},\overline{\sigma}_{xy}) = \frac{(\sigma_{xx},\sigma_{yy},\sigma_{xy})}{q_oS^2} \nonumber \\
(\overline{\sigma}_{xz},\overline{\sigma}_{yz}) = \frac{(\sigma_{xz},\sigma_{yz})}{q_oS}
\end{eqnarray}

for the applied mechanical load and by:
\begin{eqnarray}
\renewcommand{\arraystretch}{2}
(\hat{u},\hat{v}) = \frac{1}{h \alpha_m T_o S}(u,v) \nonumber \\
\hat{w} = \frac{1}{h \alpha_m T_o S^2}(w) \nonumber \\
(\hat{\sigma}_{xx},\hat{\sigma}_{yy},\hat{\sigma}_{xy}) = \frac{(\sigma_{xx},\sigma_{yy},\sigma_{xy})}{E_m \alpha_m T_o} \nonumber \\
(\hat{\sigma}_{xz},\hat{\sigma}_{yz}) = \frac{(\sigma_{xz},\sigma_{yz})}{E_m \alpha_m T_o}
\end{eqnarray}

for the applied thermal load, where the reference value is taken as $E_o =$ 1 GPa, $S = a/h$ and $E_m, \alpha_m$ are the Young's modulus and the co-efficient of thermal expansion corresponding to the metallic phase. Based on a progressive mesh-refinement, an 8 $\times$ 8 mesh is found to be adequate to model the full FGM plate for the present analysis. Before proceeding to the detailed analysis of the static response of the FGM sandwich plate for the applied mechanical and thermal loads, the present formulation is validated considering the problems for which three-dimensional elasticity solutions are available. In Table~\ref{fgmvalidation} the converged displacements and stresses obtained for the Al/SiC functionally graded square plate under mechanical {\bf loading} is compared with the three-dimensional elasticity solutions~\cite{velbatra2002}. The effective material properties are based on the Mori-Tanaka homogenization scheme. The material properties for Al are: $E_m=$ 70 GPa, $\nu_m=$ 0.3, $\alpha_m = $ 23.4 $\times 10^{-6}$ /K and $\kappa_m = $ 233 W/mK. For SiC, the material properties are: $E_m=$ 427 GPa, $\nu_m=$ 0.17, $\alpha_m = $ 4.3 $\times 10^{-6}$ /K and $\kappa_m = $ 65 W/mK. It can be seen that the results from the present formulation are in good agreement with the 3D elasticity solutions. 

Next, the numerical study is carried out for Type A FGM sandwich plate under mechanical/thermal loads for two gradient indices, $n$ (= 0,0.5) and for different thickness of the plate. Numerical results are tabulated in Tables \ref{TypeA_mechGrad05} - \ref{TypeA_mechGrad5} for an applied mechanical load and in Tables \ref{TypeA_thermGrad05} - \ref{TypeA_thermGrad5} for an applied thermal load. It is inferred that with increasing gradient index, the non-dimensional displacements and stresses increases, whereas, the displacements and the stresses decreases with increasing core thickness. This is attributed to the change in the flexural stiffness of the plates due to the increase in metallic and ceramic volume fractions. Furthermore, the displacements/stresses predicted at the neutral surface using HSDT11 and HSDT13 are somewhat different from those of HSDT9 and FSDT5. However, the latter models cannot predict the through thickness variation of displacements/stresses.

\begin{figure}[htpb]
\centering
\subfigure[$\overline{u}$]{\includegraphics[scale=0.40]{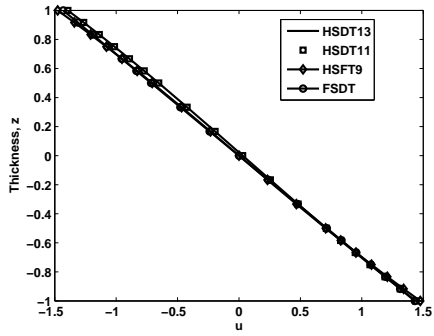}}
\subfigure[$\overline{w}$]{\includegraphics[scale=0.40]{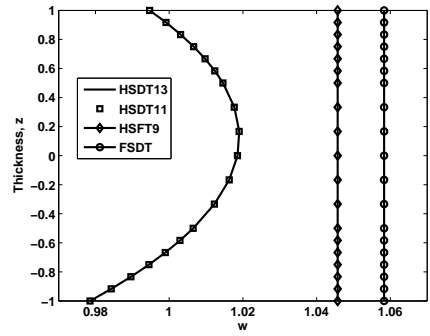}}
\subfigure[$\overline{\sigma}_{xx}$]{\includegraphics[scale=0.40]{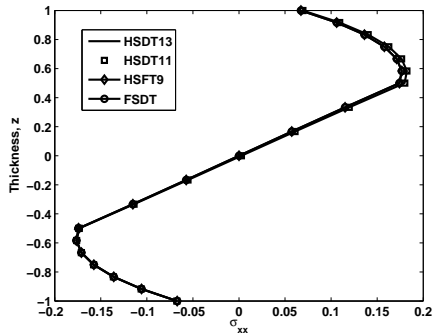}}
\subfigure[$\overline{\sigma}_{xz}$]{\includegraphics[scale=0.40]{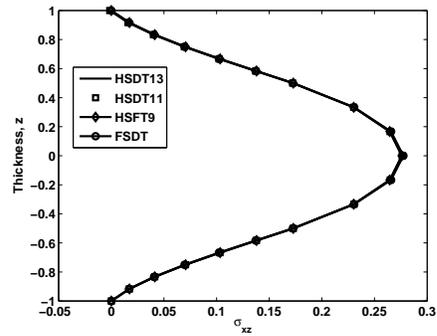}}
\caption{Displacements and stresses through the thickness for the square plates with simply supported edges for 1-2-1 Type A FGM plate with gradient index, $n =$ 1, $a/h=$ 5, subjected to a uniformly distributed mechanical loading.}
\label{fig:mefgm121_dispstress}
\end{figure}

\begin{figure}[htpb]
\centering
\subfigure[$\hat{u}$]{\includegraphics[scale=0.42]{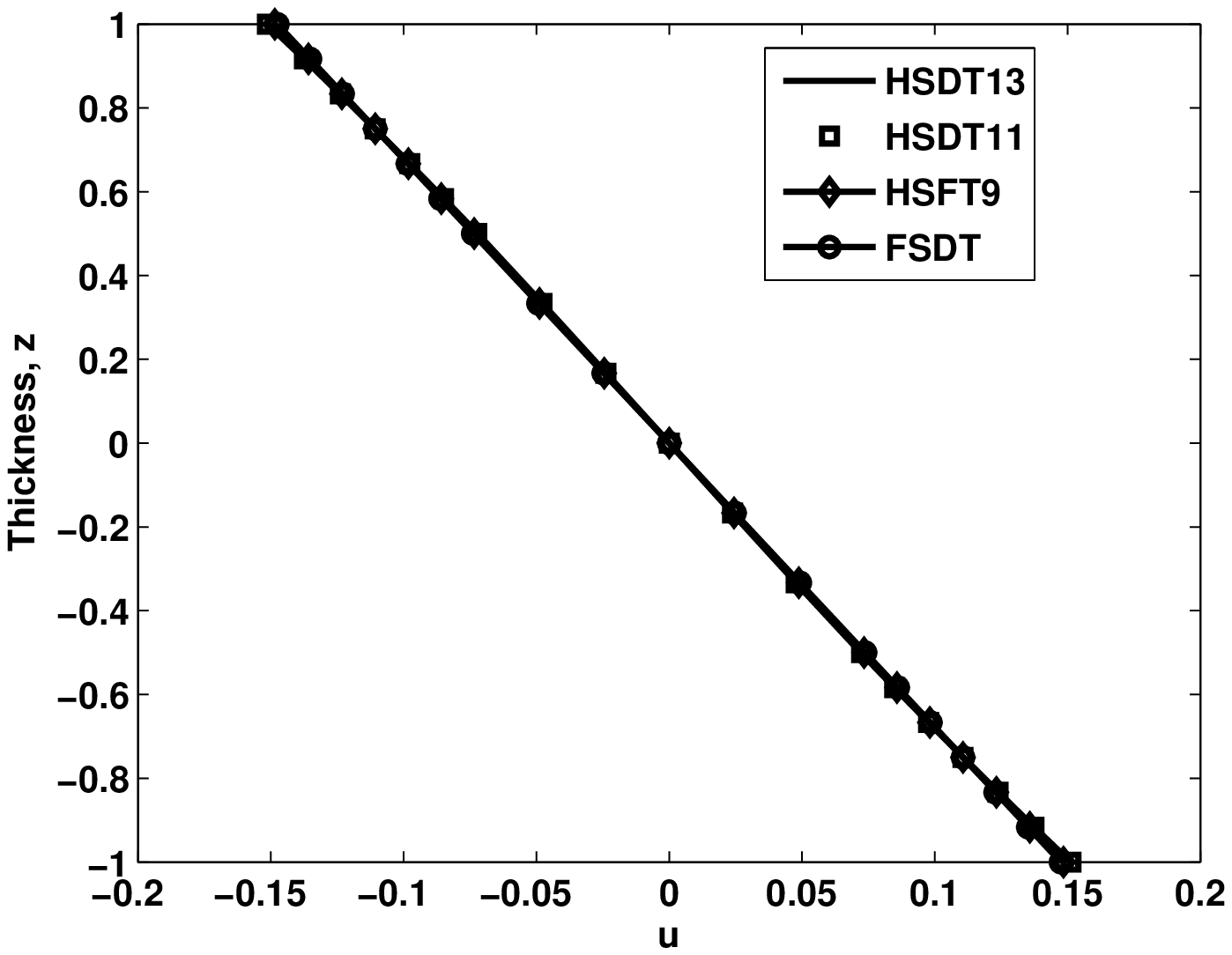}}
\subfigure[$\hat{w}$]{\includegraphics[scale=0.42]{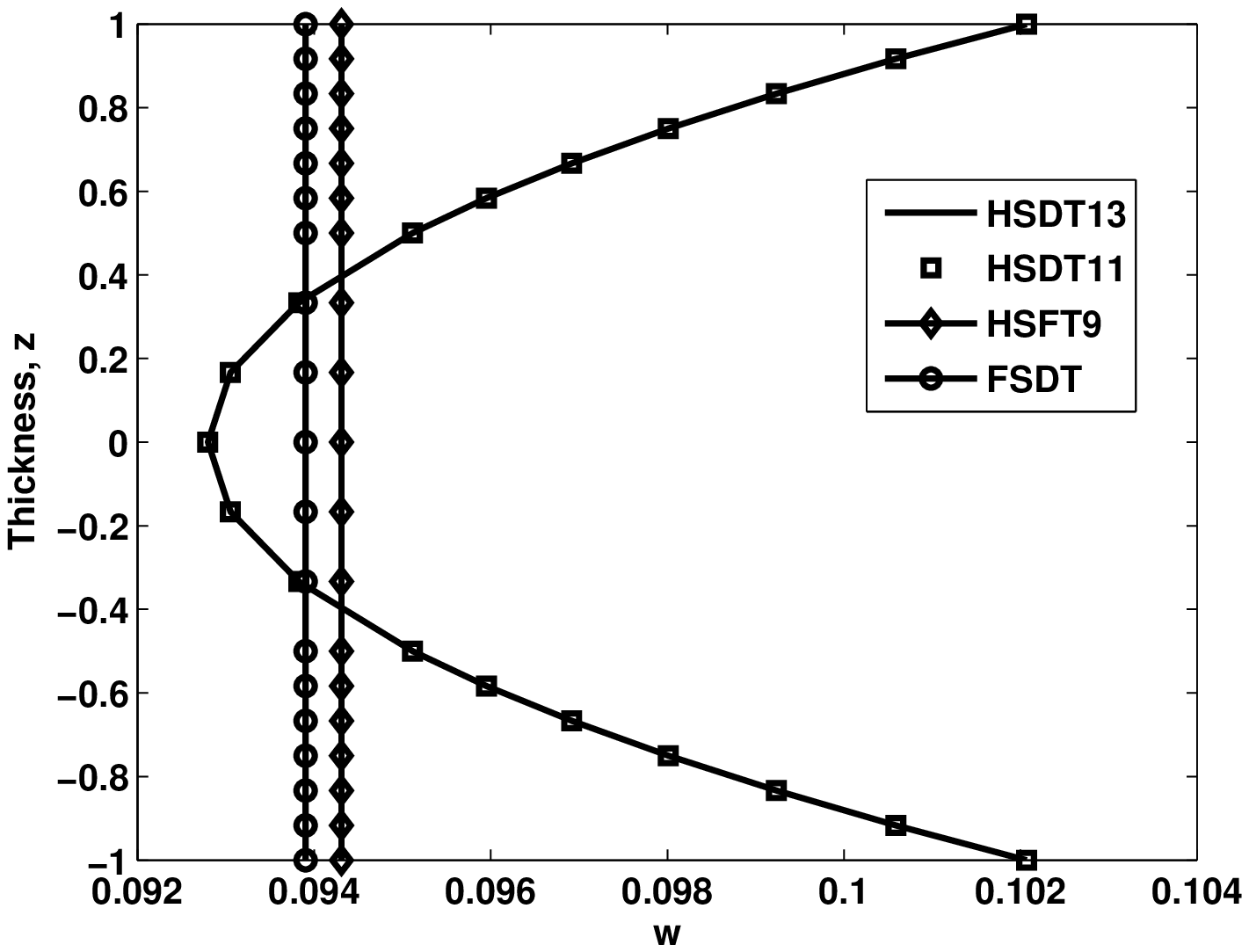}}
\subfigure[$\hat{\sigma}_{xx}$]{\includegraphics[scale=0.42]{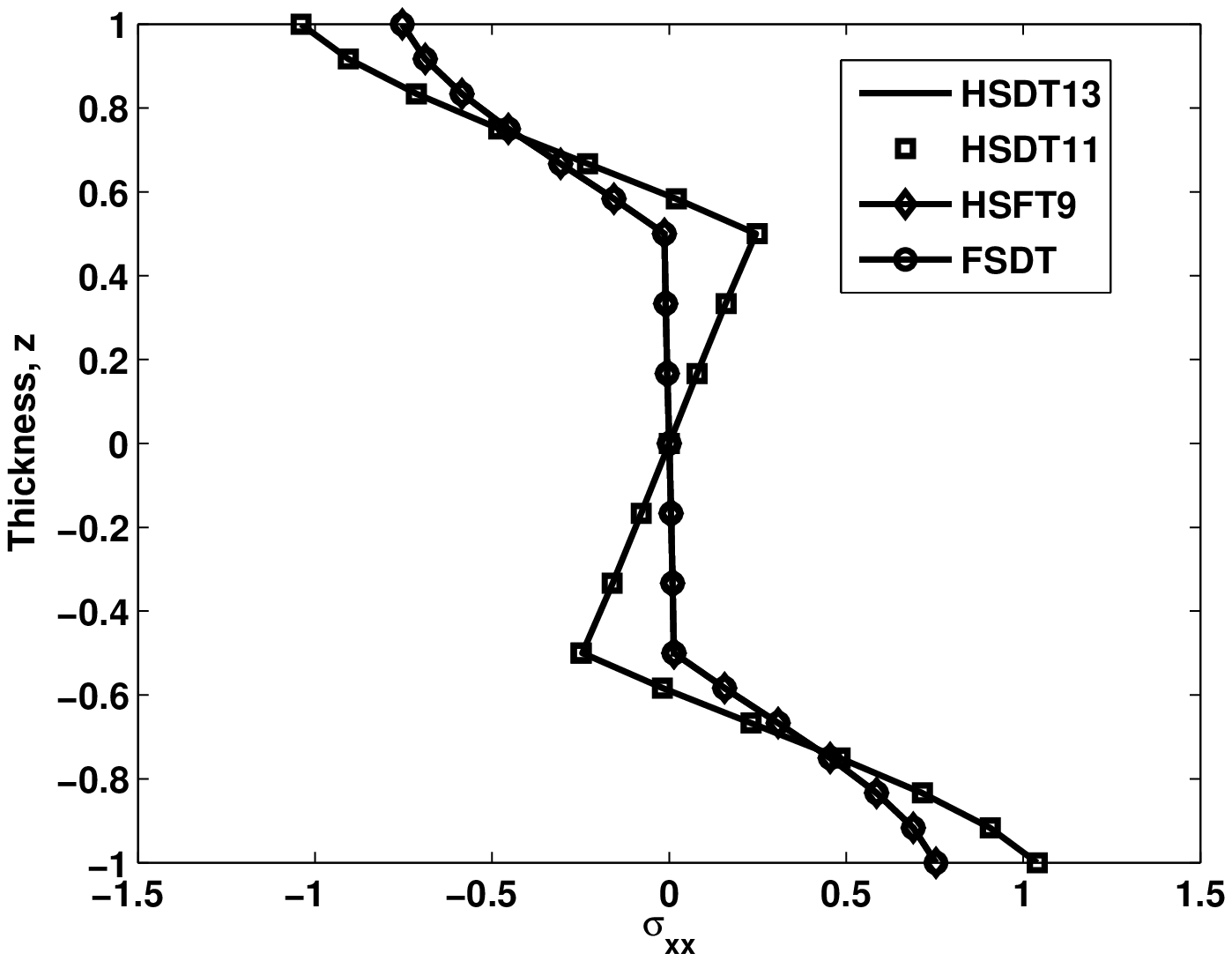}}
\subfigure[$\hat{\sigma}_{xz}$]{\includegraphics[scale=0.42]{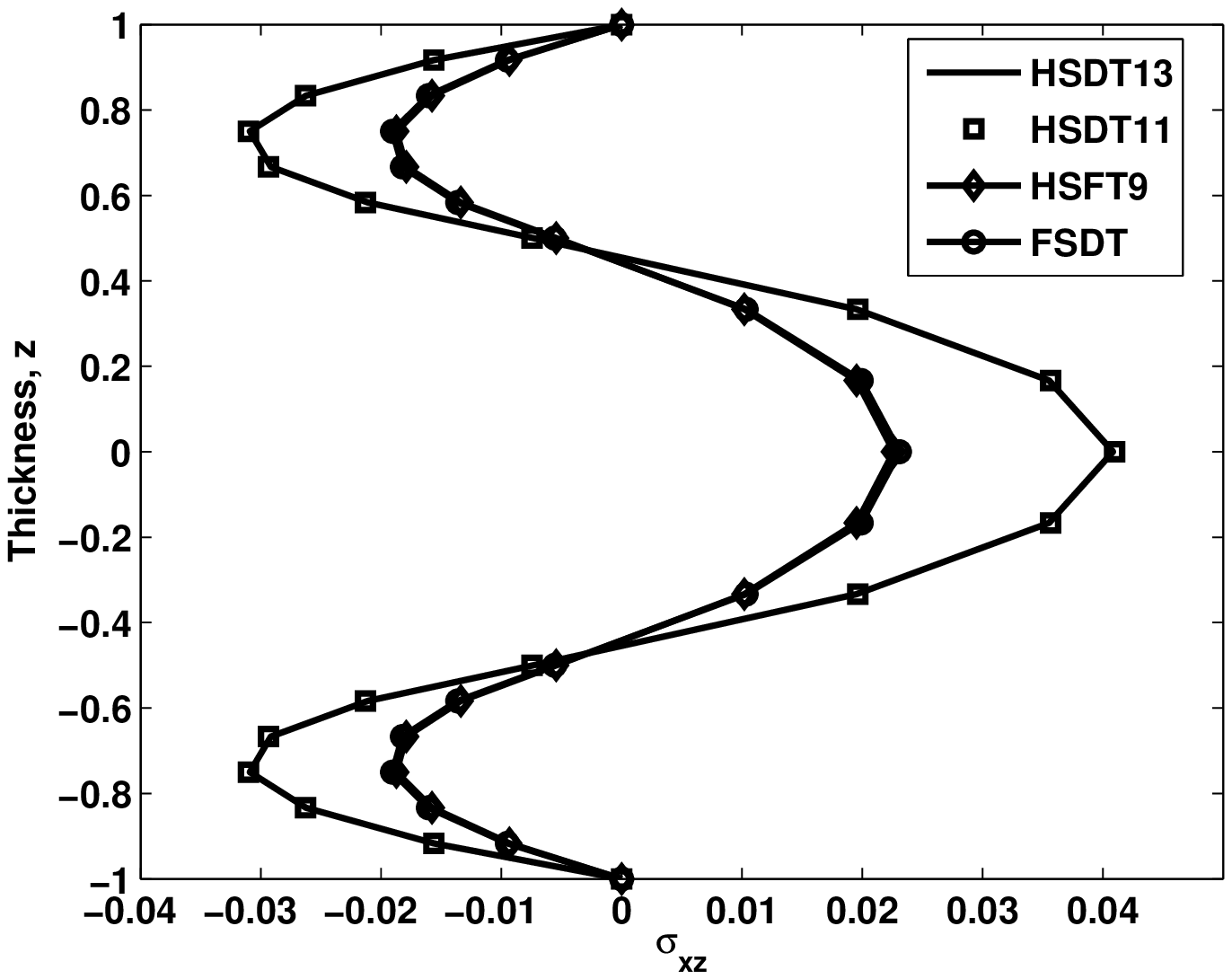}}
\caption{Displacements and stresses through the thickness for the square plates with simply supported edges for 1-2-1 Type A FGM plate with gradient index, $n =$ 1, $a/h=$ 5, subjected to a thermal loading.}
\label{fig:thfgm121_dispstress}
\end{figure}

Through the thickness variation of displacements and stresses for 1-2-1 Type A FGM sandwich plate is shown in \frefs{fig:mefgm121_dispstress} - \ref{fig:thfgm121_dispstress} for different models proposed in this study for mechanical and thermal load cases, respectively. For mechanical loading case, the higher-order (HSDT13 and HSDT11) and lower-order models (HSDT9 and FSDT5) yield almost similar results. The main difference among the higher-order elements is in accounting for the slope discontinuity in the in-plane displacements through the thickness, whereas, it is the inclusion of variation up to cubic in the in-plane displacements among lower-order elements. Since the material properties are smoothly varying through the thickness in FGM plates, there is no variation in the evaluated results seen among the higher-order and lower-order models. However, the performance of these types of elements are significantly different for thermal case (see \fref{fig:thfgm121_dispstress}), in particular, stress variation and this can be attributed to the variation of the coefficient of thermal expansion in the thickness direction. This effect is captured due to the presence of linear/quadratic terms in lateral deflection. 


\begin{table}[htpb]
\caption{The effective material properties are based on the Mori-Tanaka homogenization scheme. The results are based on the HSDT11 element and the gradient index,  $n=$ 1.}
\renewcommand{\arraystretch}{1.25}
\centering
\subtable[Comparison of displacements, stresses for the Al/SiC functionally graded square plate when subjected to a mechanical load. The nondimensionalized physical quantities are: $\overline{u} = 100E_m u/(q_ohS^3), \overline{w} = 100E_m w/(q_ohS^4), (\overline{\sigma}_{xx},\overline{\sigma}_{xy}) = 10(\sigma_{xx}, \sigma_{xy})/(q_oS^3), \overline{\sigma}_{xz}=10\sigma_{xz}/(q_oS)$]{
\centering
\begin{tabular}{clrrrrr}
\hline
$a/h$ & Mesh & \multicolumn{5}{c}{Varaible} \\
\cline{3-7}
& & $\overline{u}$ & $\overline{w}$ & $\overline{\sigma}_{xx}$ & $\overline{\sigma}_{xy}$ & $\overline{\sigma}_{xz}$\\
\hline
\multirow{5}{*}{5} & 4 $\times$ 4 & -2.9305 & 2.5407 & 2.7905 & -1.6545 & 2.2783 \\
& 6 $\times$ 6 & -2.9153 & 2.5512 & 2.7697 & -1.5960 & 2.2985\\
& 8 $\times$ 8 & -2.9129 & 2.5535 & 2.7549 & -1.5783 & 2.3016\\
& 16 $\times$ 16 & -2.9129 & 2.5535 & 2.7549 & -1.5783 & 2.3016\\
& Ref.~\cite{velbatra2002} & -2.9129 & 2.5559 & 2.7562 & -1.5600 & 2.3100\\
\cline{2-7}
\multirow{5}{*}{40} & 4 $\times$ 4 & -2.8972 & 2.1169 & 2.6564 & -1.6470 & 2.3253 \\
& 6 $\times$ 6 & -2.8972 & 2.1148 & 2.5762 & -1.5880 & 2.3287\\
& 8 $\times$ 8 & -2.8967 & 2.1152 & 2.5494 & -1.5704 & 2.3286\\
& 16 $\times$ 16 & -2.8967 & 2.1152 & 2.5494 & -1.5704 & 2.3286\\
& Ref.~\cite{velbatra2002} & -2.8984 & 2.1163 & 2.6093 & -1.5522 & 2.3281\\
\hline
\end{tabular}
} \\

\subtable[Comparison of displacements, stresses for the Al/SiC functionally graded square plate when subjected to a thermal load. The nondimensionalized physical quantities are: $\hat{u} = 100u/(h\alpha_m T_o S), \hat{w} = 100w/(h\alpha_mT_oS^2), (\hat{\sigma}_{xx},\hat{\sigma}_{xy},\hat{\sigma}_{xz}) = 10(\sigma_{xx}, \sigma_{xy},\sigma_{xz})/(E_m\alpha_mT_o)$]{
\centering
\begin{tabular}{clrrrrr}
\hline
$a/h$ & Mesh & \multicolumn{5}{c}{Varaible} \\
\cline{3-7}
& & $\hat{u}$ & $\hat{w}$ & $\hat{\sigma}_{xx}$ & $\hat{\sigma}_{xy}$ & $\hat{\sigma}_{xz}$\\
\hline
\multirow{4}{*}{5} & 4 $\times$ 4 & -1.3176 & 4.4726 & -1.2966 & -6.6678 & 2.0430\\
& 6 $\times$ 6 & -1.3232 & 4.4946 & -1.6720 & -7.2500 & 2.3848\\
& 8 $\times$ 8 & -1.3240 & 4.4499 & -1.7445 & -7.1786 & 2.5352\\
& 16 $\times$ 16 & -1.3240 & 4.4499 & -1.7445 & -7.1786 & 2.5352\\
\cline{2-7}
\multirow{4}{*}{40} & 4 $\times$ 4 & -1.2802 & 3.3660 & -1.6954 & -7.2616 & 2.3157\\
& 6 $\times$ 6 & -1.2868 & 3.3795 & -2.0216 & -7.0540 & 2.6674\\
& 8 $\times$ 8 & -1.2878 & 3.3818 & -2.1340 & -6.9832 & 2.8544 \\
& 16 $\times$ 16 & -1.2878 & 3.3818 & -2.1340 & -6.9832 & 2.8544 \\
\hline
\end{tabular}
}
\label{fgmvalidation}
\end{table}

\begin{landscape}

\begin{table}[htpb]
\centering
\renewcommand\arraystretch{1.5}
\caption{Deflections and stresses for a simply supported square FGM sandwich plates with homogeneous core (Type A), with gradient index, $n=$ 0.5, subjected to a sinusoidally distributed load (mechanical load).}
\begin{tabular}{clrrrrrrrrrrr}
\hline
$a/h$  & Element & \multicolumn{5}{c}{1-1-1} & & \multicolumn{5}{c}{1-2-1} \\
\cline{3-7} \cline{9-13} \\
&Type & $\overline{u}$ & $\overline{w}$ & $\overline{\sigma}_{xx}$ & $\overline{\sigma}_{xy}$ & $\overline{\sigma}_{xz}$ & & $\overline{u}$ & $\overline{w}$ & $\overline{\sigma}_{xx}$ & $\overline{\sigma}_{xy}$ & $\overline{\sigma}_{xz}$ \\
\hline
\multirow{4}{*}{5}& HSDT13  &  0.01827  &  0.01257  &  -0.05962  &  0.03131  &  0.26308  &&  0.01677  &  0.01158  &  -0.05469  &  0.02872  &  0.26010  \\
&HSDT11  &  0.01827  &  0.01257  &  -0.05964  &  0.03131  &  0.26321  &&  0.01677  &  0.01158  &  -0.05471  &  0.02872  &  0.26010  \\
&HSDT9  &  0.01867  &  0.01342  &  -0.05942  &  0.03199  &  0.26284  &&  0.01712  &  0.01236  &  -0.05447  &  0.02933  &  0.25980  \\
&FSDT5  &  0.01813  &  0.01358  &  -0.05769  &  0.03106  &  0.26449  &&  0.01656  &  0.01247  &  -0.05269  &  0.02837  &  0.26158  \\
\cline{2-13}                                            
\multirow{4}{*}{10}&HSDT13 & 0.01818 & 0.01181  &  -0.05822  &  0.03114  &  0.26438  &&  0.01662  &  0.01081  &  -0.05323  &  0.02847  &  0.26150  \\
&HSDT11 & 0.01817 & 0.01181 & -0.05823 &  0.03114  &  0.26447  &&  0.01662  &  0.01081  &  -0.05324  &  0.02847  &  0.26150  \\
&HSDT9 & 0.01827 & 0.01201 & -0.05813  & 0.03130 & 0.26413 && 0.01670 & 0.01099 & -0.05313 & 0.02861 & 0.26119 \\
&FSDT5 & 0.01813 & 0.01205 & -0.05769 & 0.03106 & 0.26450 && 0.01656 & 0.01102 & -0.05269 & 0.02837 & 0.26159 \\
\cline{2-13}                      
\multirow{4}{*}{100}& HSDT13 & 0.01813 & 0.01154 & -0.05771 & 0.03107 & 0.26548 && 0.01655 & 0.01054 & -0.05270 & 0.02838 & 0.26255 \\
&HSDT11 & 0.01813 & 0.01154 & -0.05771 & 0.03107 & 0.26548 && 0.01655 & 0.01054 & -0.05270 & 0.02838 & 0.26255 \\
&HSDT9 & 0.01813 & 0.01154 & -0.05769 & 0.03107 & 0.26515 && 0.01655 & 0.01054 & -0.05268 & 0.02838 & 0.26222 \\
&FSDT5 & 0.01813 & 0.01155 & -0.05768 & 0.03107 & 0.26514 && 0.01655 & 0.01054 & -0.05268 & 0.02837 & 0.26221 \\
\hline
\end{tabular}
\label{TypeA_mechGrad05}
\end{table}

\begin{table}[htpb]
\centering
\renewcommand\arraystretch{1.5}
\caption{Deflections and stresses for a simply supported square FGM sandwich plates with homogeneous core (Type A), with gradient index, $n=$ 5, subjected to a sinusoidally distributed load (mechanical load).}
\begin{tabular}{clrrrrrrrrrrr}
\hline
$a/h$  & Element & \multicolumn{5}{c}{1-1-1} & & \multicolumn{5}{c}{1-2-1} \\
\cline{3-7} \cline{9-13} \\
& Type & $\overline{u}$ & $\overline{w}$ & $\overline{\sigma}_{xx}$ & $\overline{\sigma}_{xy}$ & $\overline{\sigma}_{xz}$  & & $\overline{u}$ & $\overline{w}$ & $\overline{\sigma}_{xx}$ & $\overline{\sigma}_{xy}$ & $\overline{\sigma}_{xz}$ \\
\hline
\multirow{4}{*}{5}&HSDT13  &  0.04232  &  0.02828  &  -0.13876  &  0.07250  &  0.31653  &&  0.03233  &  0.02151  &  -0.10626  &  0.05538  &  0.31370  \\
&HSDT11  &  0.04230  &  0.02827  &  -0.13871  &  0.07247  &  0.31626  &&  0.03232  &  0.02151  &  -0.10626  &  0.05537  &  0.31360  \\
&HSDT9  &  0.04362  &  0.03024  &  -0.13883  &  0.07474  &  0.31516  &&  0.03332  &  0.02303  &  -0.10602  &  0.05708  &  0.31277  \\
&FSDT5  &  0.04353  &  0.03077  &  -0.13852  &  0.07457  &  0.31534  &&  0.03277  &  0.02340  &  -0.10430  &  0.05615  &  0.31376  \\
\cline{2-13}                                  
\multirow{4}{*}{10}&HSDT13  &  0.04323  &  0.02785  &  -0.13860  &  0.07407  &  0.31601  &&  0.03267  &  0.02102  &  -0.10483  &  0.05598  &  0.31413  \\
&HSDT11 &  0.04323  &  0.02784  &  -0.13860  &  0.07406  &  0.31599  &&  0.03267  &  0.02102  &  -0.10483  &  0.05597  &  0.31409  \\
&HSDT9  &  0.04355  &  0.02834  &  -0.13859  &  0.07462  &  0.31534  &&  0.03291  &  0.02140  &  -0.10473  &  0.05639  &  0.31356  \\
&FSDT5  &  0.04353  &  0.02847  &  -0.13852  &  0.07458  &  0.31537  &&  0.03277  &  0.02150  &  -0.10430  &  0.05615  &  0.31378  \\
\cline{2-13}
\multirow{4}{*}{100}&HSDT13  &  0.04351  &  0.02771  &  -0.13853  &  0.07460  &  0.31667  &&  0.03277  &  0.02086  &  -0.10432  &  0.05617  &  0.31504  \\
&HSDT11  &  0.04351  &  0.02771  &  -0.13853  &  0.07460  &  0.31667  &&  0.03277  &  0.02086  &  -0.10432  &  0.05617  &  0.31504  \\
&HSDT9  &  0.04352  &  0.02771  &  -0.13849  &  0.07460  &  0.31626  &&  0.03277  &  0.02087  &  -0.10428  &  0.05618  &  0.31464  \\
&FSDT5  &  0.04352  &  0.02771  &  -0.13849  &  0.07460  &  0.31623  &&  0.03277  &  0.02087  &  -0.10428  &  0.05617  &  0.31461  \\

\hline
\end{tabular}
\label{TypeA_mechGrad5}
\end{table}

\end{landscape}


\begin{landscape}

\begin{table}[htpb]
\centering
\renewcommand\arraystretch{1.25}
\caption{Deflections and stresses for a simply supported square FGM sandwich plates with homogeneous core (Type A), with gradient index, $n=$ 0.5, subjected to a sinusoidally distributed load (thermal loading).}
\begin{tabular}{clrrrrrrrrrrr}
\hline
$a/h$  & Element& \multicolumn{5}{c}{1-1-1} & & \multicolumn{5}{c}{1-2-1} \\
\cline{3-7} \cline{9-13} \\
& Type& $\hat{u}$ & $\hat{w}$ & $\hat{\sigma}_{xx}$ & $\hat{\sigma}_{xy}$ & $\hat{\sigma}_{xz}$ & & $\hat{u}$ & $\hat{w}$ & $\hat{\sigma}_{xx}$ & $\hat{\sigma}_{xy}$ & $\hat{\sigma}_{xz}$\\
\hline
\multirow{4}{*}{5}&HSDT13	& 0.14143 & 0.09551 & 1.13358 & 0.34643 & -0.01628 && 0.13275 & 0.08966 & 1.21775 & 0.32518 & -0.00813 \\
&HSDT11 & 0.14148 & 0.09551 & 1.13318 & 0.34657 & -0.01637 && 0.13280 & 0.08965 & 1.21740 & 0.32529 & -0.00826 \\
&HSDT9 & 0.13894 & 0.08825 & 0.79687 & 0.34030 & -0.01026 && 0.13038 & 0.08282 & 0.83580 & 0.31933 & -0.00732 \\
&FSDT5 & 0.13821 & 0.08799 & 0.80020 & 0.33845 & -0.01046 && 0.12962 & 0.08252 & 0.83924 & 0.31742 & -0.00755 \\
\cline{2-13}                      
\multirow{4}{*}{10}&HSDT13 & 0.13902 & 0.08987 & 1.14995 & 0.34052 & 0.02687 && 0.13041 & 0.08431 & 1.23367 & 0.31943 & 0.01808 \\
&HSDT11 & 0.13903 & 0.08987 & 1.14982 & 0.34055 & 0.02692 && 0.13042 & 0.08431 & 1.23355 & 0.31946 & 0.01811 \\
&HSDT9 & 0.13839 & 0.08806 & 0.79938 & 0.33895 & 0.01389 && 0.12981 & 0.08260 & 0.83839 & 0.31793 & 0.00891 \\
&FSDT5 & 0.13821 & 0.08799 & 0.80020 & 0.33846 & 0.01397 && 0.12962 & 0.08252 & 0.83924 & 0.31742 & 0.00896 \\
\cline{2-13}                      
\multirow{4}{*}{100}&HSDT13 & 0.13821 & 0.08801 & 1.15545 & 0.33851 & 0.00275 && 0.12962 & 0.08254 & 1.23900 & 0.31747 & 0.00185 \\
&HSDT11 & 0.13821 & 0.08801 & 1.15545 & 0.33851 & 0.00275 && 0.12962 & 0.08254 & 1.23900 & 0.31747 & 0.00185 \\
&HSDT9 & 0.13821 & 0.08799 & 0.80022 & 0.33850 & 0.00142 && 0.12962 & 0.08253 & 0.83927 & 0.31747 & 0.00092 \\
&FSDT5 & 0.13820 & 0.08799 & 0.80023 & 0.33849 & 0.00142 && 0.12962 & 0.08252 & 0.83927 & 0.31746 & 0.00092 \\
\hline
\end{tabular}
\label{TypeA_thermGrad05}
\end{table}

\begin{table}[htpb]
\centering
\renewcommand\arraystretch{1.25}
\caption{Deflections and stresses for a simply supported square FGM sandwich plates with homogeneous core (Type A), with gradient index, $n=$ 5, subjected to a sinusoidally distributed load (thermal loading).}
\begin{tabular}{clrrrrrrrrrrr}
\hline
$a/h$  & Element & \multicolumn{5}{c}{1-1-1} & & \multicolumn{5}{c}{1-2-1} \\
\cline{3-7} \cline{9-13} \\
& Type& $\hat{u}$ & $\hat{w}$ & $\hat{\sigma}_{xx}$ & $\hat{\sigma}_{xy}$ & $\hat{\sigma}_{xz}$ & & $\hat{u}$ & $\hat{w}$ & $\hat{\sigma}_{xx}$ & $\hat{\sigma}_{xy}$ & $\hat{\sigma}_{xz}$\\
\hline
\multirow{4}{*}{5}&HSDT13 & 0.18232 & 0.12238 & 0.73973 & 0.44652 & -0.02475 && 0.15938 & 0.10641 & 0.96614 & 0.39045 & -0.05416 \\
&HSDT11 & 0.18240 & 0.12247 & 0.73915 & 0.44679 & -0.02497 && 0.15942 & 0.10644 & 0.96585 & 0.39059 & -0.05439 \\
&HSDT9 & 0.17878 & 0.11317 & 0.61575 & 0.43787 & -0.01539 && 0.15555 & 0.09815 & 0.72136 & 0.38102 & -0.03237 \\
&FSDT5 & 0.17771 & 0.11314 & 0.62060 & 0.43519 & -0.01570 && 0.15382 & 0.09793 & 0.72923 & 0.37668 & -0.03299 \\
\cline{2-13}                   
\multirow{4}{*}{10}&HSDT13 & 0.17888 & 0.11544 & 0.76316 & 0.43812 & 0.02160 && 0.15522 & 0.10004 & 0.99415 & 0.38025 & 0.01437 \\
&HSDT11 & 0.17889 & 0.11547 & 0.76301 & 0.43819 & 0.02166 && 0.15523 & 0.10005 & 0.99408 & 0.38029 & 0.01442 \\
&HSDT9 & 0.17798 & 0.11315 & 0.61940 & 0.43590 & 0.01152 && 0.15425 & 0.09798 & 0.72728 & 0.37783 & 0.00732 \\
&FSDT5 & 0.17771 & 0.11314 & 0.62060 & 0.43519 & 0.01159 && 0.15382 & 0.09793 & 0.72924 & 0.37668 & 0.00736 \\
\cline{2-13}                      
\multirow{4}{*}{100}&HSDT13 & 0.17771 & 0.11317 & 0.77116 & 0.43527 & 0.00221 && 0.15382 & 0.09795 & 1.00359 & 0.37676 & 0.00147 \\
&HSDT11 & 0.17771 & 0.11317 & 0.77116 & 0.43527 & 0.00221 && 0.15382 & 0.09795 & 1.00359 & 0.37676 & 0.00147 \\
&HSDT9 & 0.17771 & 0.11314 & 0.62064 & 0.43525 & 0.00117 && 0.15381 & 0.09793 & 0.72926 & 0.37674 & 0.00075 \\
&FSDT5 & 0.17770 & 0.11314 & 0.62064 & 0.43524 & 0.00117 && 0.15381 & 0.09793 & 0.72927 & 0.37672 & 0.00075 \\
\hline
\end{tabular}
\label{TypeA_thermGrad5}
\end{table}

\end{landscape}


\subsection{Free flexural vibrations}
Next, the free vibration characteristics of FGM sandwich plate is numerically studied. In all cases, we present the non-dimensionalized free flexural frequency defined as:

\begin{equation}
\Omega = \frac{\omega a^2}{h} \sqrt{ \frac{ \rho_o}{E_o}}
\label{eqn:vibnondim}
\end{equation}

where $\omega$ is the natural frequency, $\rho_o =$ 1 kg/m$^3$ and $E_o=$1 GPa. Before proceeding with the detailed study on the effect of gradient index and the type of sandwich FGM plates on the natural frequencies, the formulation developed is validated against available three-dimensional elasticity solutions~\cite{Li2008}. Based on a progressive refinement, a 8$\times$8 mesh is found to be sufficient to model the sandwich FGM plate. Table \ref{TypeA_Conver} gives a comparison of the first six computed frequencies for a simply supported square Type A FGM sandwich plate. It can be seen that the numerical results from the present formulation are found to be in very good agreement with the existing solutions.

\begin{table}[htpb]
\centering
\renewcommand\arraystretch{1.5}
\caption{Convergence of flexural vibration frequency parameters $\Omega$ of square 2-1-2 FGM plates of Type A.}
\begin{tabular}{cclrrrrrr}
\hline
$n$ & $a/h$ & Mesh & \multicolumn{6}{c}{Frequency}\\
\cline{4-9}
 &  & & $\Omega_1$ & $\Omega_2$ & $\Omega_3$ & $\Omega_4$ & $\Omega_5$ & $\Omega_6$ \\
\hline
\multirow{9}{*}{1} & \multirow{4}{*}{5} & 4 $\times$ 4 & 1.2297 & 2.6874 & 2.6876 & 2.8149 & 2.8271 & 4.1312 \\
& & 6 $\times$6 & 1.2294 & 2.6869 & 2.6869 & 2.8032 & 3.2551 & 4.1553 \\
& & 8 $\times$8 & 1.2293 & 2.6868 & 2.6868 & 2.8009 & 2.8345 & 4.1568 \\
& & 16 $\times$16 & 1.2293 & 2.6868 & 2.6868 & 2.8009 & 2.8345 & 4.1568 \\
\cline{3-9}
& \multirow{4}{*}{10} & 4 $\times$ 4 & 1.3025 & 3.1822 & 3.1822 & 4.9167 & 6.2782 & 7.7331 \\
& & 6 $\times$6 & 1.3020 & 3.1640 & 3.1640 & 4.9209 & 6.0938 & 7.6898 \\
& & 8 $\times$8 & 1.3019 & 3.1606 & 3.1606 & 4.9188 & 6.0586 & 7.7675 \\
& & 16 $\times$16 & 1.3019 & 3.1606 & 3.1606 & 4.9188 & 6.0586 & 7.7675 \\
& & Li \textit{et al,}~\cite{Li2008} & 1.3018 & 3.1588 & 3.1588 & 4.9166 & 6.0405& --  \\
\cline{2-9}
\multirow{9}{*}{10} & \multirow{4}{*}{5} & 4 $\times$ 4 & 0.8958 & 2.0729 & 2.0729 & 2.2067 & 2.2225 & 3.0703 \\
& & 6$\times$6 & 0.8955 & 2.0636 & 2.0637 & 2.2063 & 2.5395 & 3.0840 \\
& & 8$\times$8 & 0.8955 & 2.0618 & 2.0619 & 2.2062 & 2.4973 & 3.0845 \\
& & 16$\times$16 & 0.8955 & 2.0618 & 2.0619 & 2.2062 & 2.4973 & 3.0845 \\
\cline{3-9}
& \multirow{4}{*}{10} & 4 $\times$ 4 & 0.9423 & 2.3109 & 2.3109 & 3.5841 & 4.5875 & 5.6835 \\
& & 6$\times$6 & 0.9419 & 2.2974 & 2.2974 & 3.5851 & 4.4488 & 5.6311 \\
& & 8$\times$8 & 0.9418 & 2.2948 & 2.2948 & 3.5832 & 4.4225 & 5.6183 \\
& & 16$\times$16 & 0.9418 & 2.2948 & 2.2948 & 3.5832 & 4.4225 & 5.6183 \\
& & Li \textit{et al,}~\cite{Li2008} & 0.9404 & 2.2862 & 2.2862 & 3.5647 & 4.3844 & -- \\
\hline
\end{tabular}
\label{TypeA_Conver}
\end{table}

\begin{landscape}

\begin{table}[htpb]
\centering
\renewcommand\arraystretch{1.5}
\caption{Fundamental frequency parameters $\Omega$ of simply supported square FGM sandwich plates with homogeneous core (Type A).}
\begin{tabular}{clrrrrcrrrcrrr}
\hline
$a/h$  & Element & \multicolumn{4}{c}{1-1-1} & & \multicolumn{3}{c}{1-2-1} & & \multicolumn{3}{c}{2-2-1}\\
\cline{3-6} \cline{8-10} \cline{12-14} 
 & Type & 0 & 0.5 & 1 & 5 & & 0.5 & 1 & 5 & & 0.5 & 1 & 5 \\
\hline
\multirow{5}{*}{5} &HSDT13 & 1.6774 & 1.4219 & 1.2778 & 0.9986 & & 1.4696 & 1.3536 & 1.1192 & & 1.4455 & 1.3144 & 1.0565\\
&HSDT11 & 1.6774 & 1.4219 & 1.2778 & 0.9988 & & 1.4696 & 1.3537 & 1.1193 & & 1.4455 & 1.3144 & 1.0566 \\
&HSDT9 & 1.6774 & 1.4152 & 1.2714 & 0.9937 & & 1.4626 & 1.3468 & 1.1131 & &1.4387 & 1.3078 & 1.0510 \\
&FSDT & 1.6689 & 1.4076 & 1.2628 & 0.9860 & & 1.4565 & 1.3398 & 1.1053 & & 1.4320 & 1.3002 & 1.0444 \\
& Li~\textit{et al,}~\cite{Li2008}& 1.6771 & 1.4218 & 1.2777 & 0.9980 & & 1.4694 & 1.3534 & 1.1190 & & 1.4454 & 1.3143 & 1.0561 \\
\cline{2-14}
\multirow{5}{*}{10}&HSDT13 & 1.8269 & 1.5214 & 1.3553 & 1.0455 & & 1.5768 & 1.4415 & 1.1757 & & 1.5494 & 1.3977 & 1.1100 \\
&HSDT11 & 1.8269 & 1.5214 & 1.3553 & 1.0456 & & 1.5768 & 1.4415 & 1.1758 & & 1.5494 & 1.3977 & 1.1100 \\
&HSDT9 & 1.8245 & 1.5193 & 1.3553 & 1.0441 & & 1.5746 & 1.4394 & 1.1740 & &1.5472 & 1.3957& 1.1084\\
&FSDT & 1.8242 & 1.5168 & 1.3506 & 1.0418 & & 1.5726 & 1.4371 & 1.1715 & & 1.5451 & 1.3932 & 1.1064 \\
& Li~\textit{et al,}~\cite{Li2008}&1.8268 & 1.5213 & 1.3552 & 1.0453 &  &1.5767 & 1.4414 & 1.1757 & & 1.5493 & 1.3976 & 1.1098 \\
\cline{2-14}
\multirow{5}{*}{100} &HSDT13 & 1.8884 & 1.5605 & 1.3852 & 1.0631 & & 1.6192 & 1.4756 & 1.1970 & & 1.5904 & 1.4300& 1.1303 \\
&HSDT11 & 1.8884 & 1.5605 & 1.3852 & 1.0631 & &1.6192 & 1.4756 & 1.1970 & &1.5904 & 1.4300	& 1.1303 \\
&HSDT9 & 1.8883 & 1.5605 & 1.3851 & 1.0631 & & 1.6192 & 1.4756 & 1.1970  & & 1.5904 & 1.4300 &1.1302 \\
&FSDT & 1.8883 & 1.5605 & 1.3851 & 1.0631 & & 1.6192 & 1.4756 & 1.1970 & & 1.5904 & 1.4299 & 1.1302 \\
& Li~\textit{et al,}~\cite{Li2008}&1.8883 & 1.5605 & 1.3851 & 1.0631 & & 1.6192 & 1.4756 & 1.1970 & & 1.5903 & 1.4299 & 1.1302 \\
\hline
\end{tabular}
\label{TypeA_FundaFreq}
\end{table}

\begin{table}[htpb]
\centering
\renewcommand\arraystretch{1.5}
\caption{Fundamental frequency parameters $\Omega$ of simply supported square FGM sandwich plates with FGM core (Type B).}
\begin{tabular}{clrrrrcrrrcrrr}
\hline
$a/h$  & Element & \multicolumn{4}{c}{1-1-1} & & \multicolumn{3}{c}{1-2-1} & & \multicolumn{3}{c}{2-2-1}\\
\cline{3-6} \cline{8-10} \cline{12-14} 
 & Type & 0 & 0.5 & 1 & 5 & & 0.5 & 1 & 5 & & 0.5 & 1 & 5 \\
\hline
\multirow{4}{*}{5} &HSDT13 & 1.0893 & 1.1511 & 1.1701 & 1.2162 & & 1.1663 & 1.1952 & 1.2712 & & 1.2031 & 1.2421 &1.3312\\
&HSDT11 & 1.1078 & 1.1512 & 1.1705 & 1.2184 & & 1.1664 & 1.1953 & 1.2718 & & 1.2034 & 1.2422 & 1.3326 \\
&HSDT9 & 1.1021 & 1.1449 & 1.1639 & 1.2113 & &1.1597 & 1.1884 & 1.2644 & & 1.1965 & 1.2350 & 1.3249\\
&FSDT & 1.1263 & 1.1503 & 1.1642 & 1.2050 & &1.1660 & 1.1880 & 1.2567 & & 1.1950 & 1.2299 & 1.3173 \\
\cline{2-14}
\multirow{4}{*}{10} & HSDT13 & 1.2087 & 1.2392 & 1.2524 & 1.2935 & & 1.2598 & 1.2806 & 1.3513 & & 1.2865 & 1.3238 & 1.4180 \\
&HSDT11 & 1.2156 & 1.2392 & 1.2525 & 1.2942 & & 1.2598 & 1.2806 & 1.3515 & & 1.2866 & 1.3238 & 1.4184\\
&HSDT9 & 1.2138 & 1.2373 & 1.2506 & 1.2921 &  & 1.2578 & 1.2785 & 1.3492 & &1.2846 & 1.3216 & 1.4161 \\
&FSDT & 1.2225 & 1.2394 & 1.2509 & 1.2903 &  &1.2601 & 1.2786 & 1.3469 & & 1.2842 & 1.3201 & 1.4136\\
\cline{2-14}
\multirow{4}{*}{100} &HSDT13 & 1.2616 & 1.2751 & 1.2854 & 1.3239 & & 1.2981 & 1.3148 & 1.3825 & &1.3198 & 1.3559 & 1.4519 \\
&HSDT11 & 1.2617 & 1.2751 & 1.2854 & 1.3239 & &1.2981 & 1.3148 & 1.3825 & &1.3198 & 1.3559 & 1.4519 \\
&HSDT9 & 1.2617 & 1.2751 & 1.2854 & 1.3239 & &1.2981 & 1.3148 & 1.3825 & & 1.3198 & 1.3559 & 1.4519 \\
&FSDT & 1.2618 & 1.2751 & 1.2854 & 1.3239 & & 1.2981 & 1.3148 & 1.3825 & & 1.3198 & 1.3559 & 1.4518 \\
\hline
\end{tabular}
\label{TypeB_FundaFreq}
\end{table}

\end{landscape}

Tables~\ref{TypeA_FundaFreq} and \ref{TypeB_FundaFreq} gives the fundamental frequencies for Type A and Type B sandwich FGM plate for three different thickness ratios and for different core thickness, respectively. It can be seen that the natural frequencies decreases with increasing gradient index for Type A and the natural frequencies decreases with decreasing gradient index for Type B FGM sandwich plate. This can be attributed to the decrease in the material rigidity. In the case of Type A, as the gradient index increases, the metallic volume fraction increases, thus decreasing the overall stiffness of the plate. In the case of Type B, decreasing the gradient index, decreases the metallic volume fraction and increases the material rigidity. 

\begin{figure}[htpb]
\centering
\subfigure[mode 1]{\includegraphics[scale=0.42]{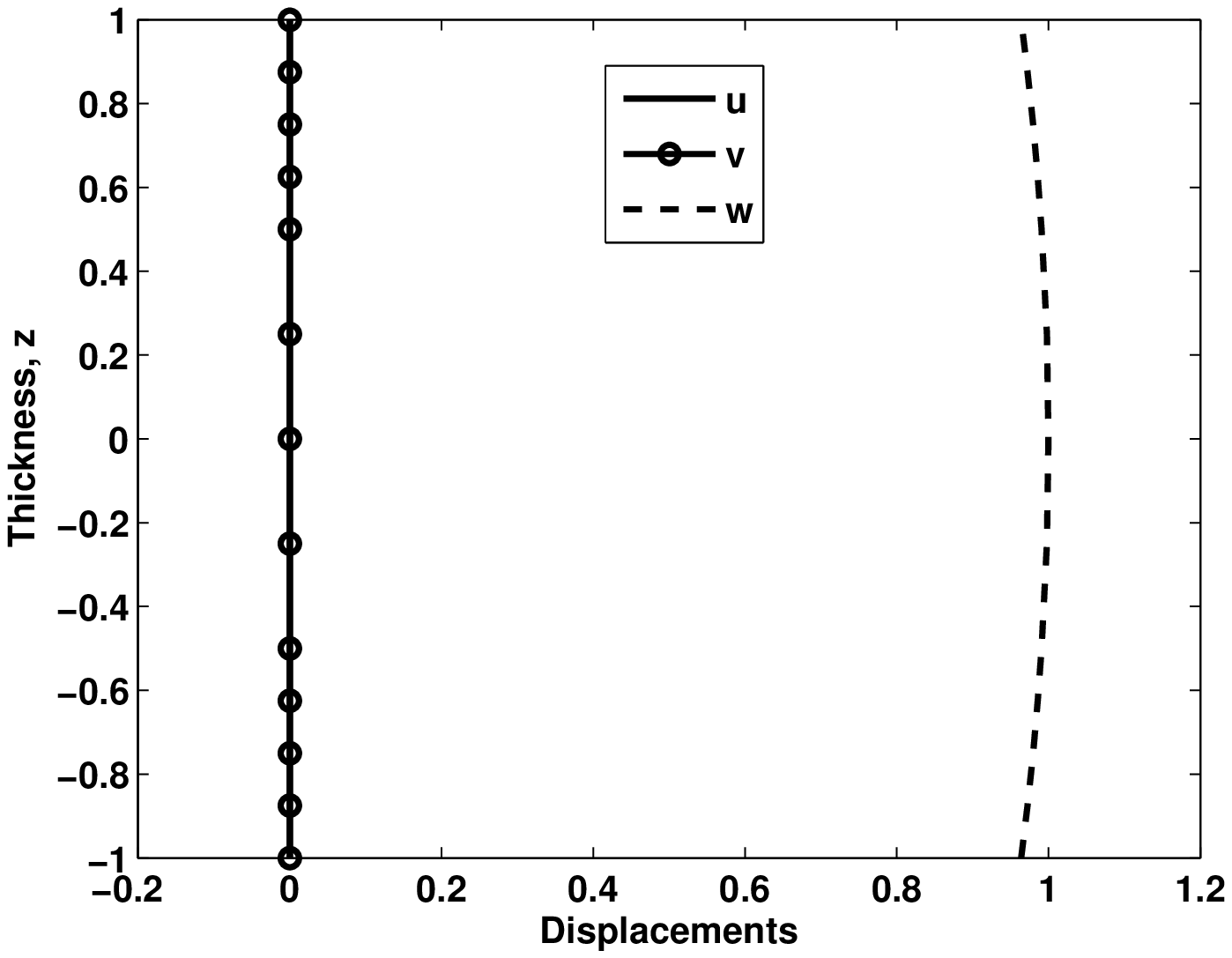}}
\subfigure[mode 2]{\includegraphics[scale=0.42]{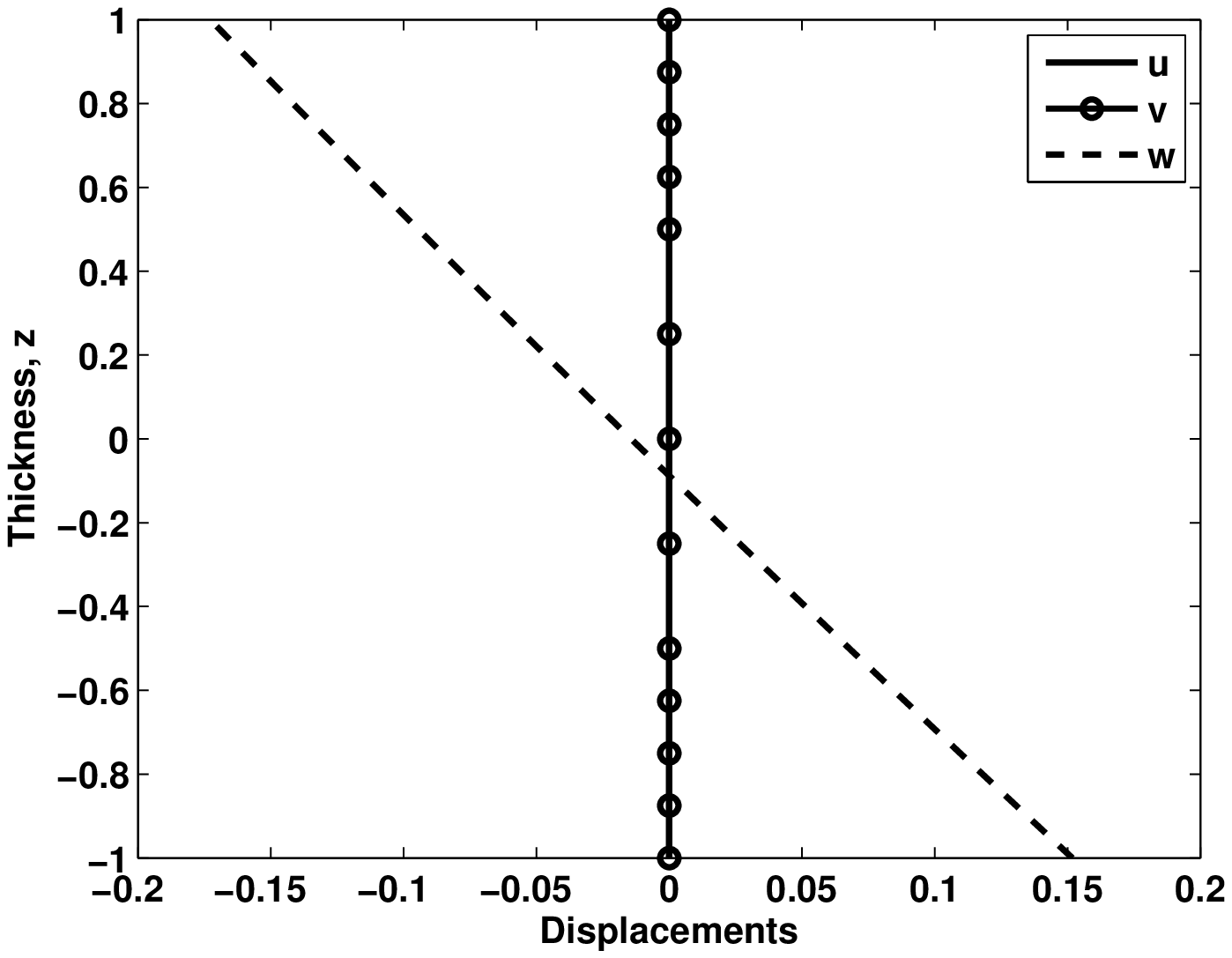}}
\subfigure[mode 3]{\includegraphics[scale=0.42]{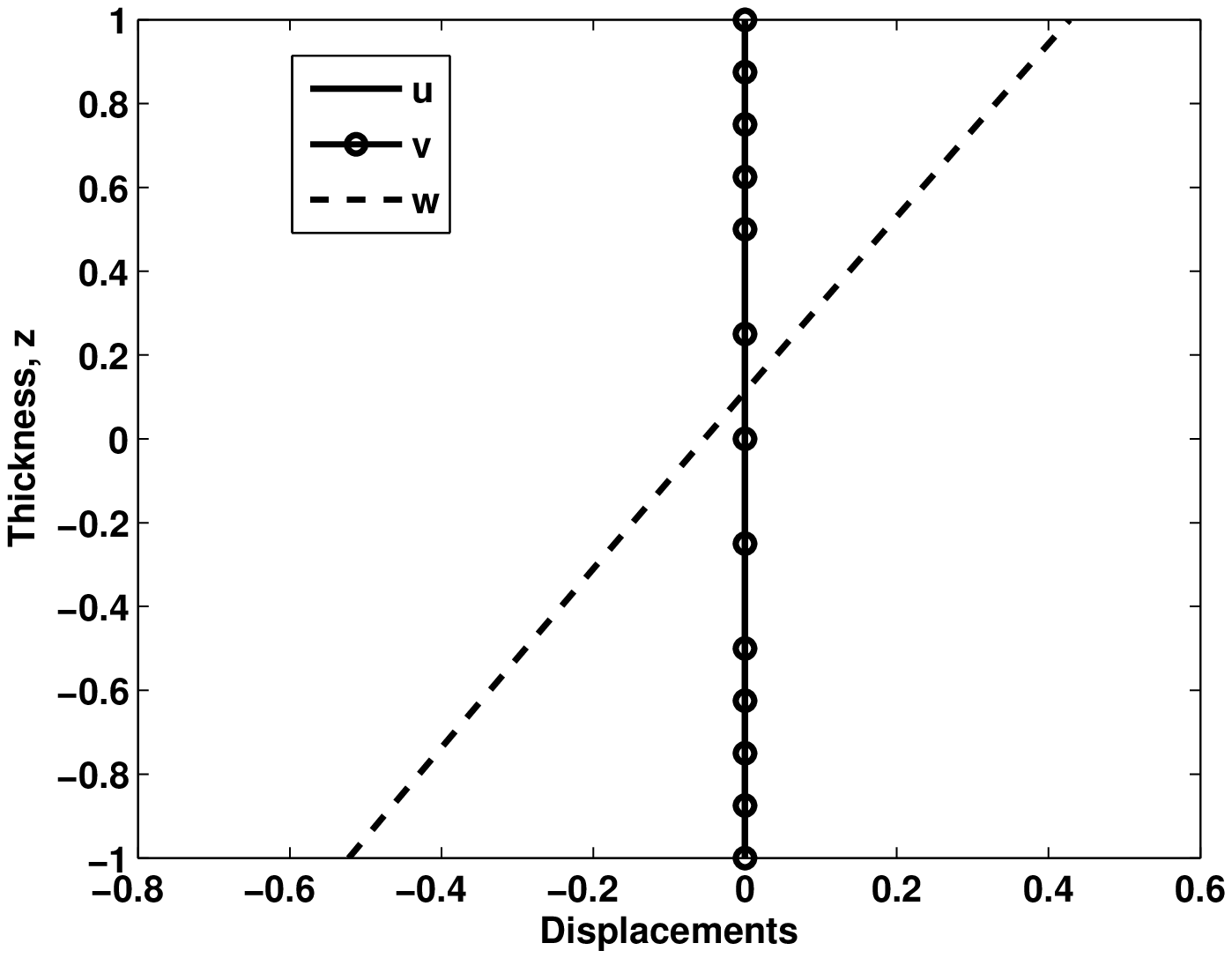}}
\subfigure[mode 4]{\includegraphics[scale=0.42]{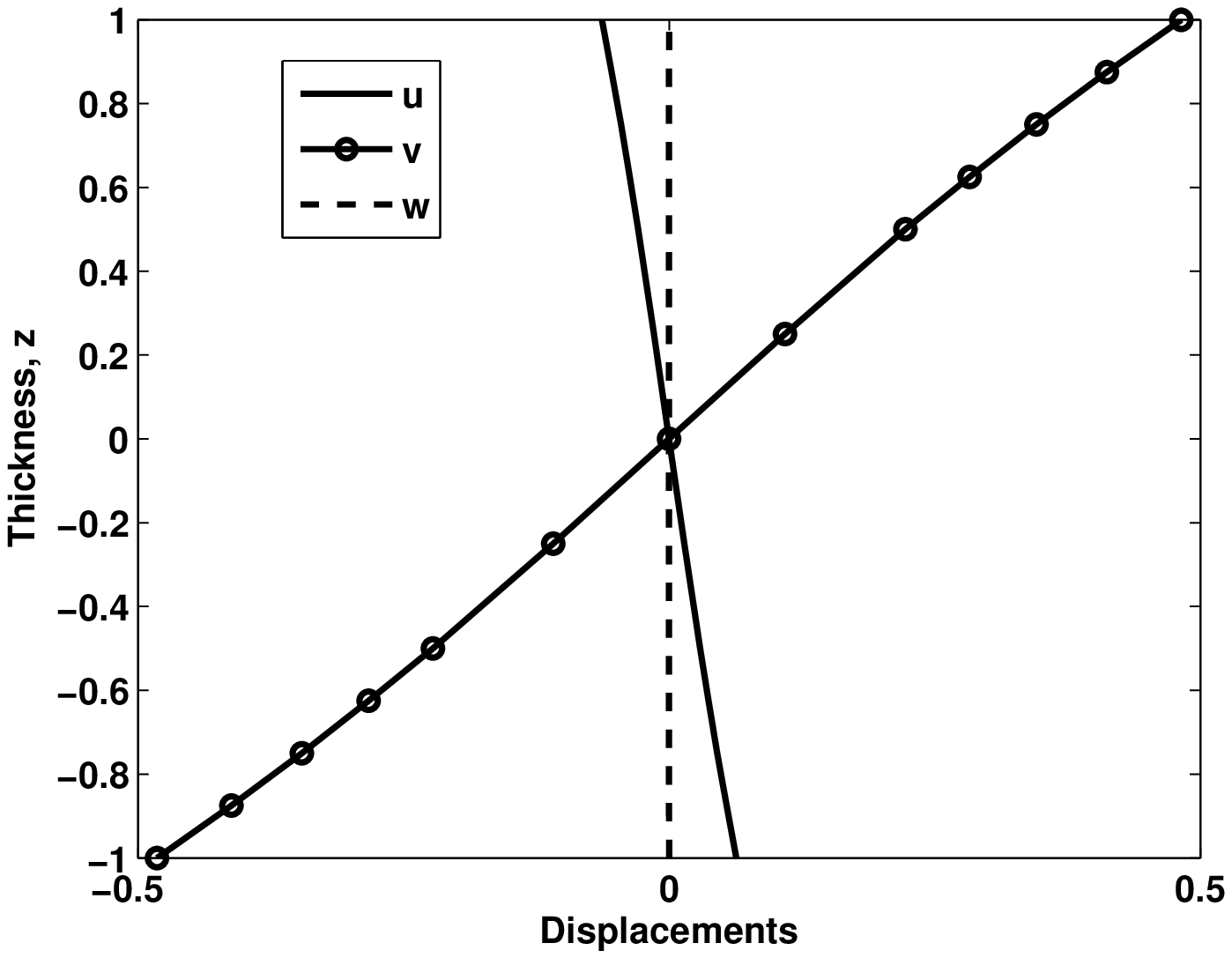}}
\subfigure[mode 5]{\includegraphics[scale=0.42]{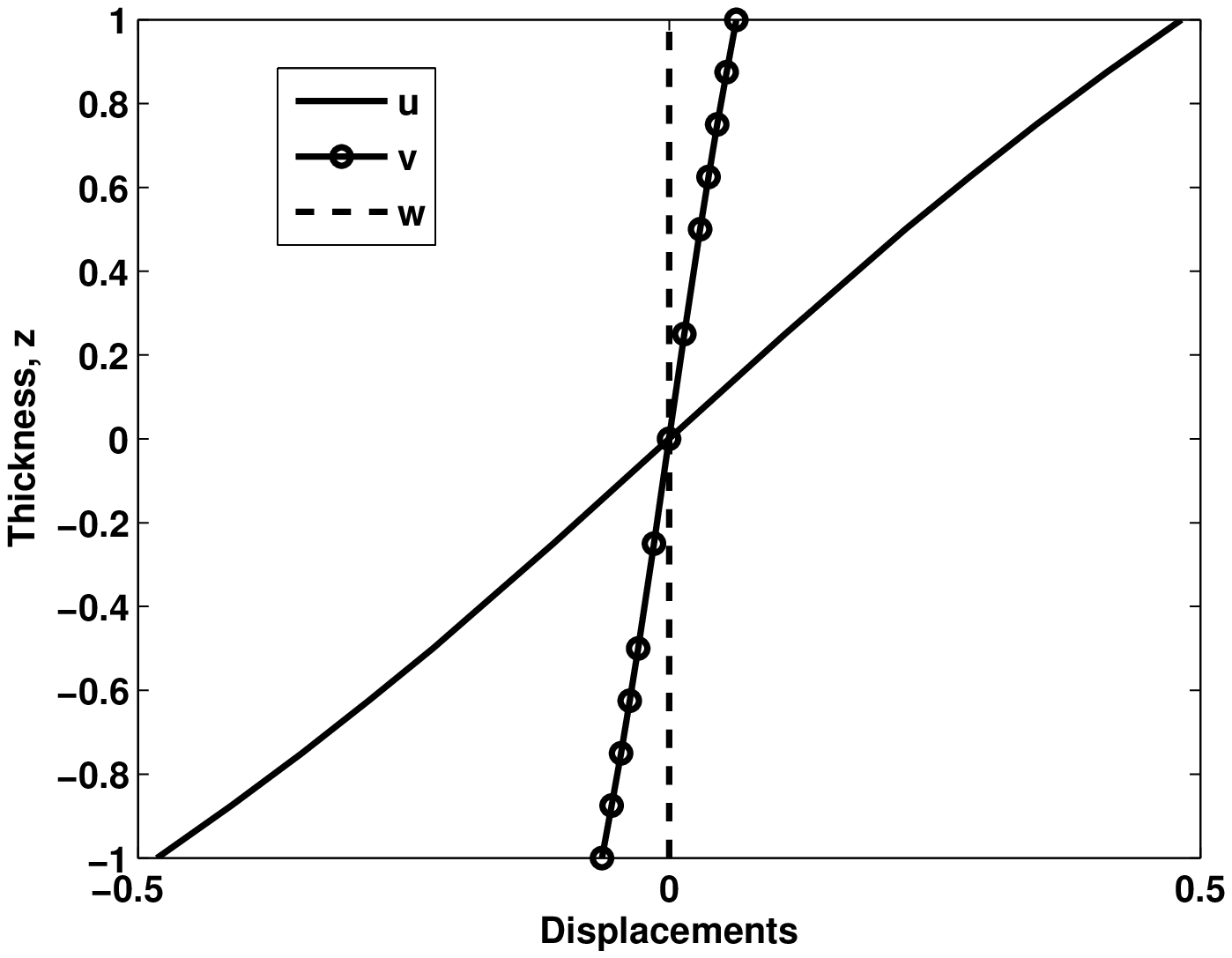}}
\subfigure[mode 6]{\includegraphics[scale=0.42]{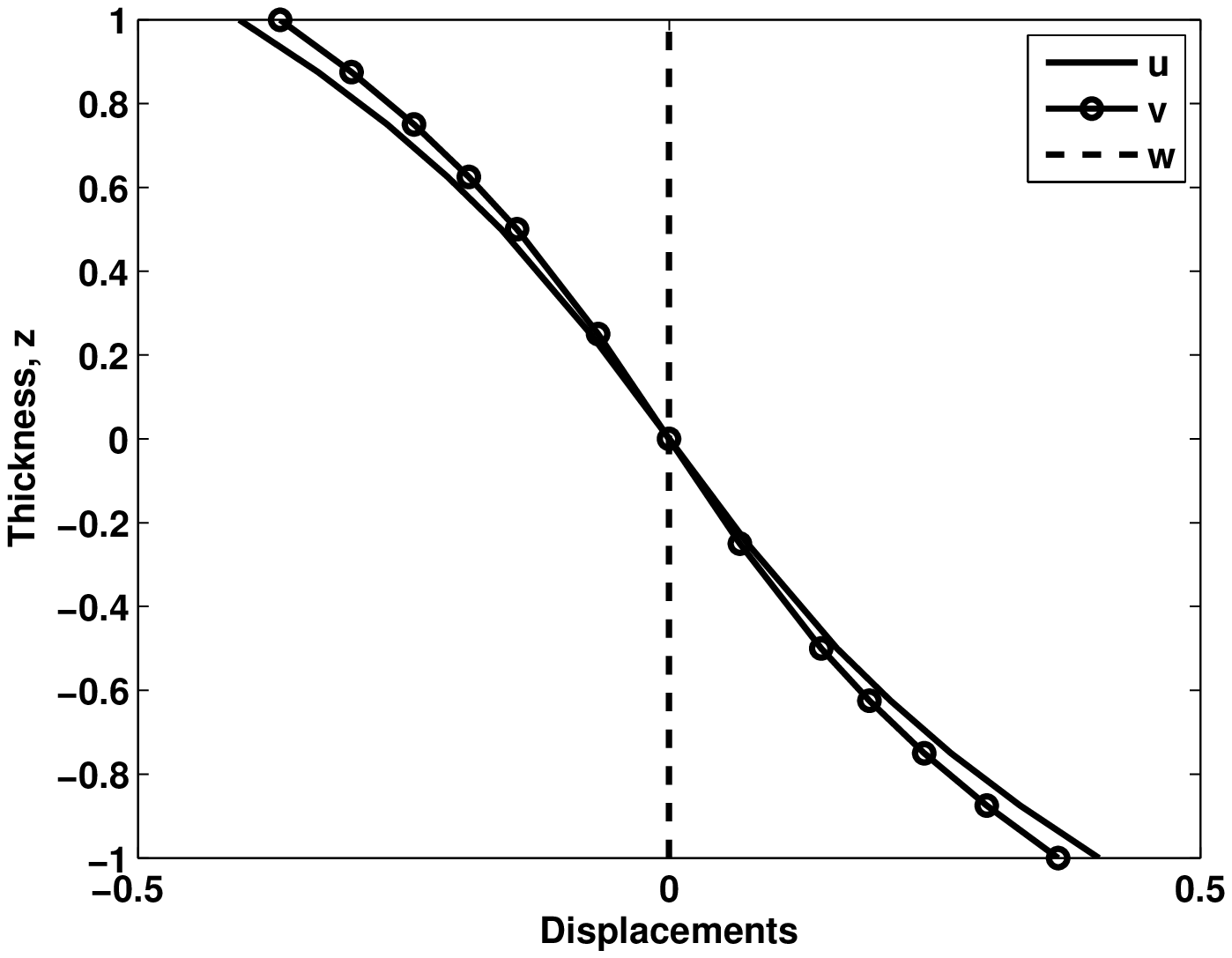}}
\caption{Deflected shapes $u_i(a/2,a/2,z),~i=x,y,z$ of the six modes for the square plates with simply supported edges for 1-2-1 Type A FGM plate with gradient index, $n =$ 1 and thickness, $a/h=$ 5.}
\label{fig:fgm121_centerA}
\end{figure}

\begin{figure}[htpb]
\centering
\subfigure[mode 1]{\includegraphics[scale=0.42]{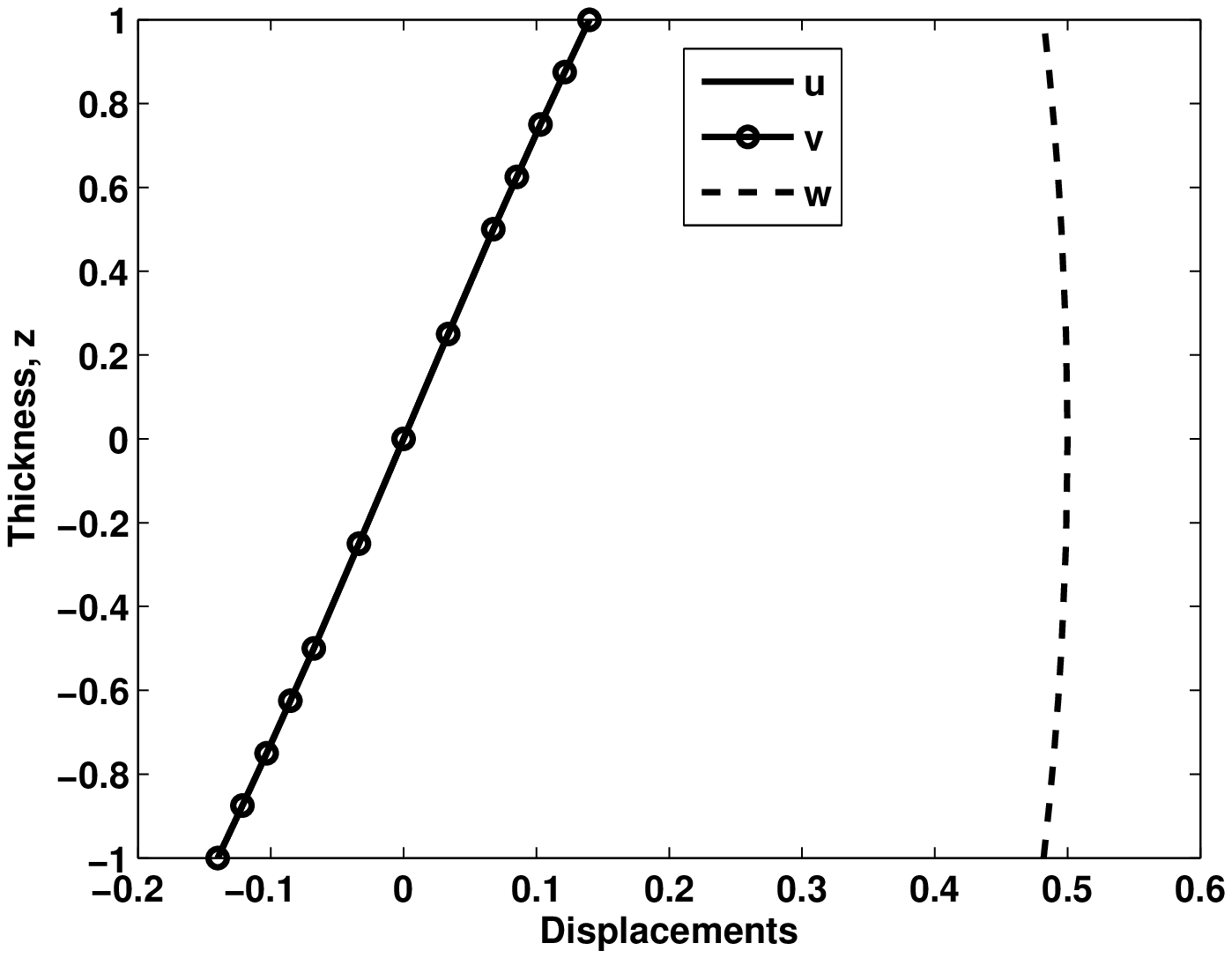}}
\subfigure[mode 2]{\includegraphics[scale=0.42]{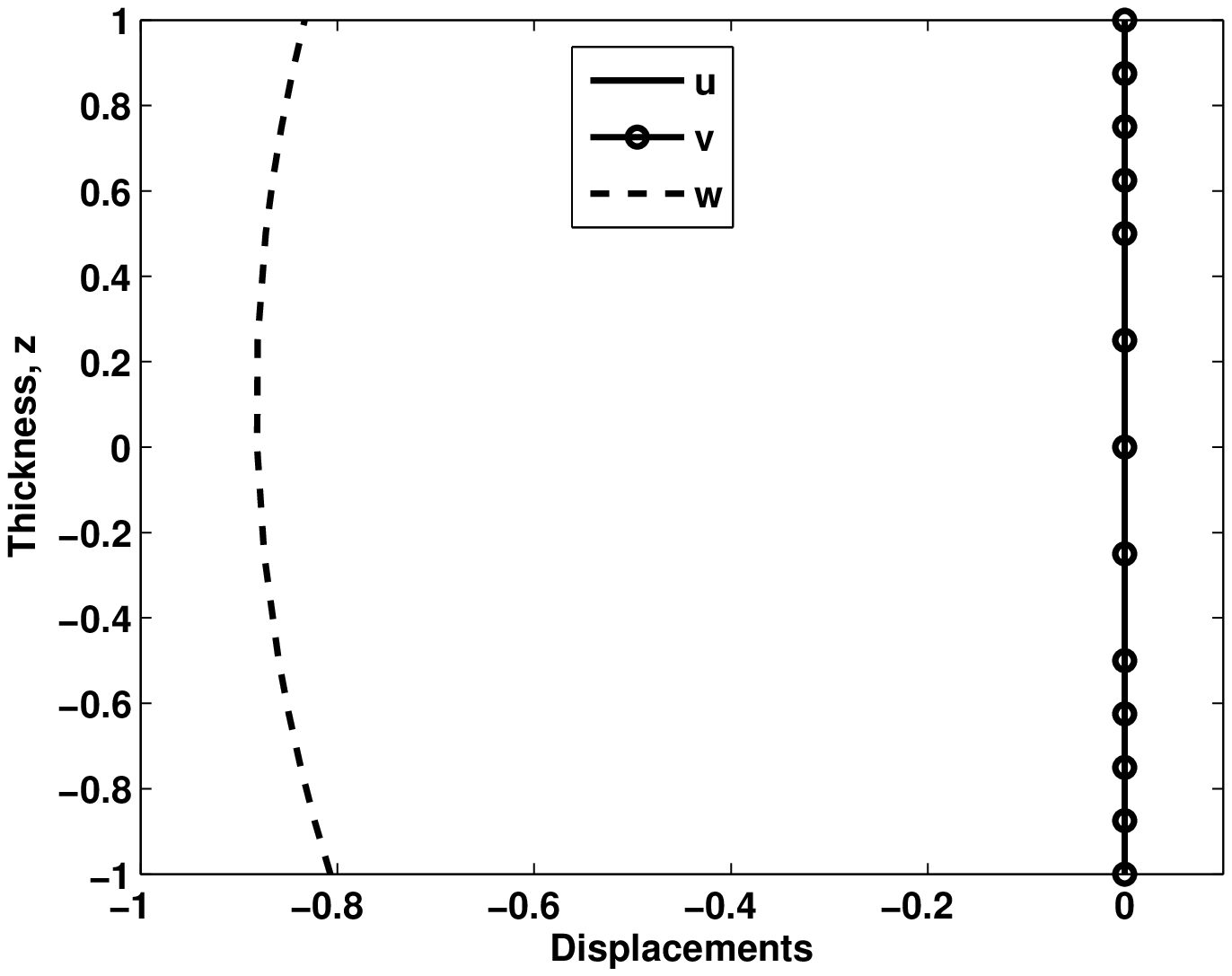}}
\subfigure[mode 3]{\includegraphics[scale=0.42]{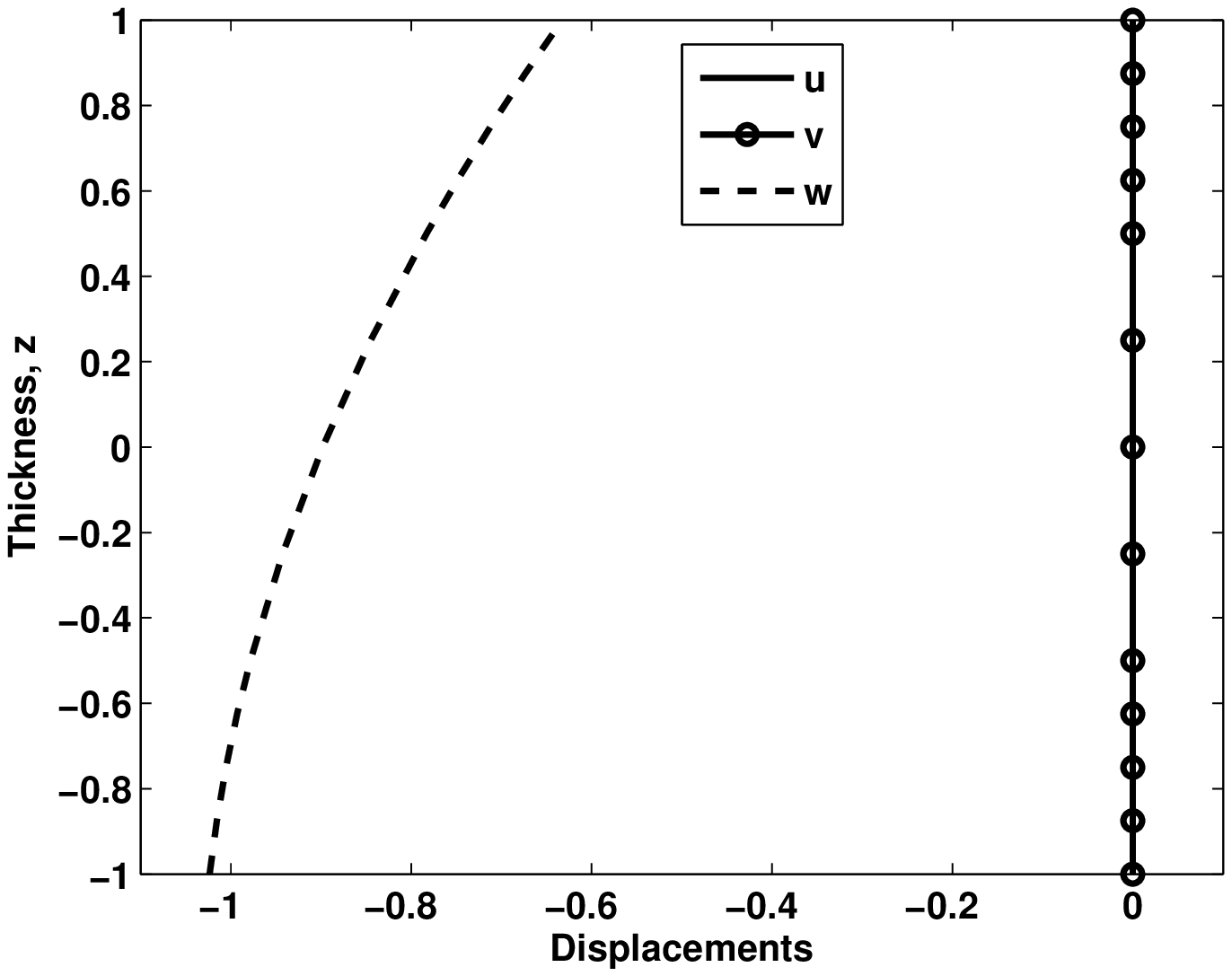}}
\subfigure[mode 4]{\includegraphics[scale=0.42]{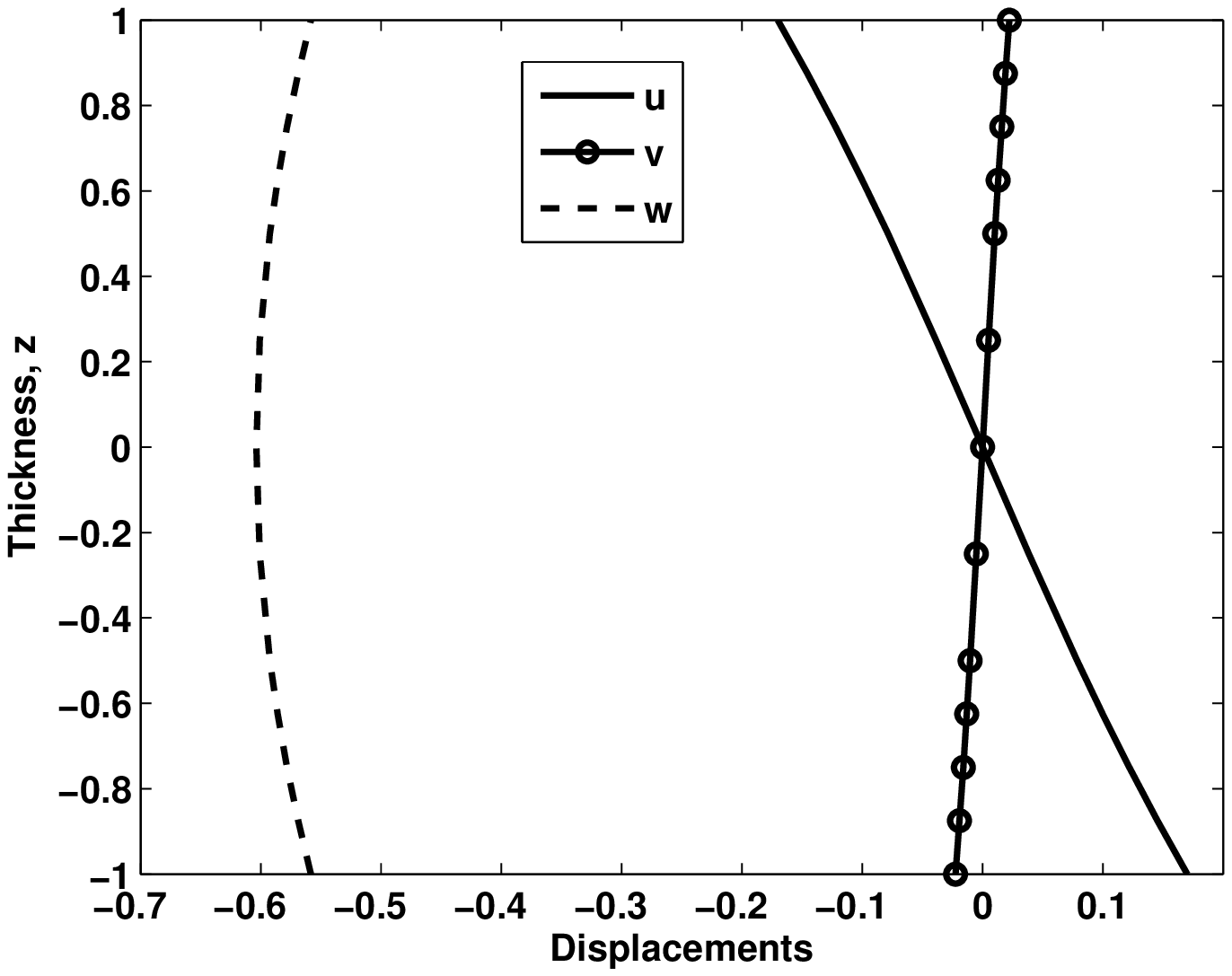}}
\subfigure[mode 5]{\includegraphics[scale=0.42]{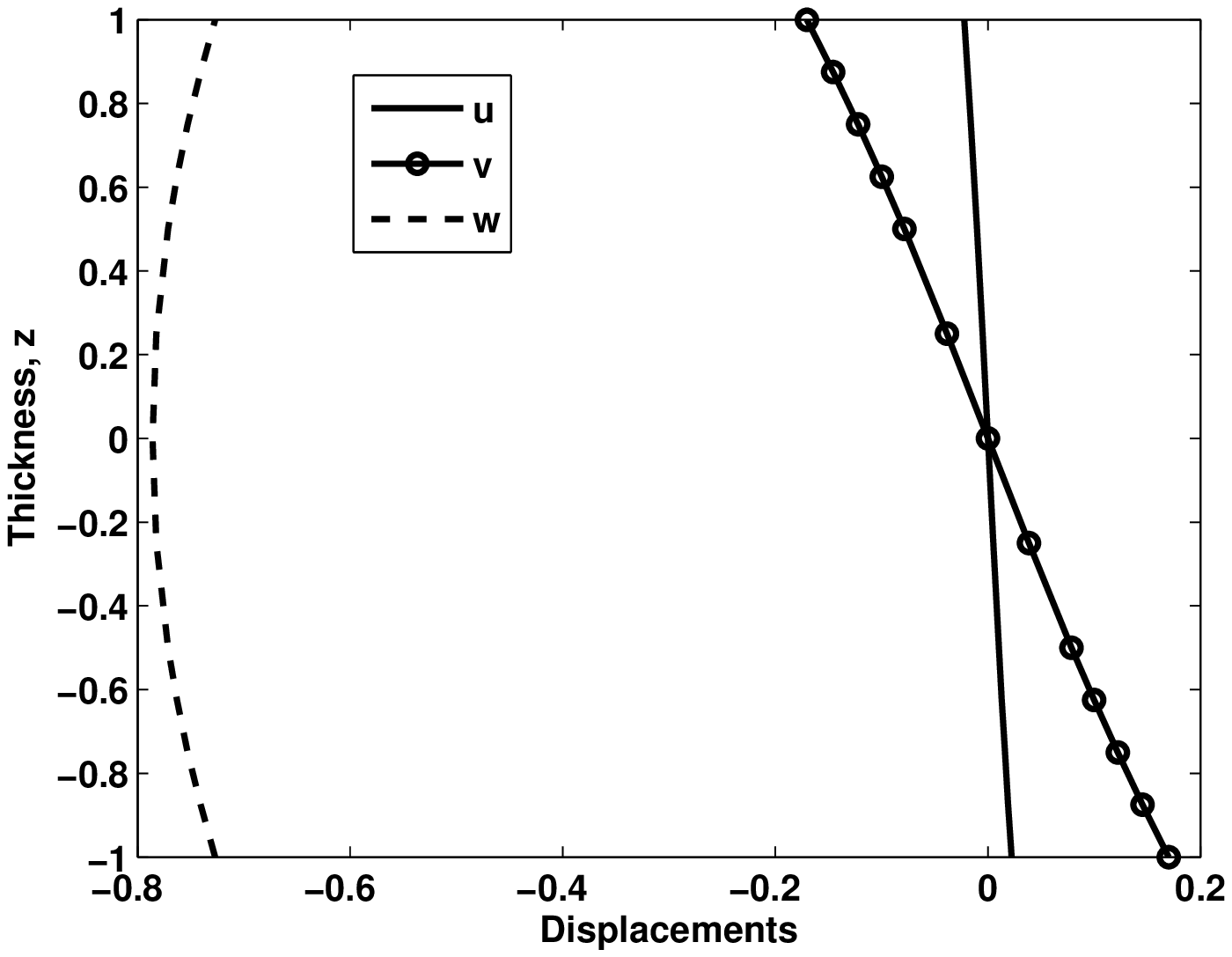}}
\subfigure[mode 6]{\includegraphics[scale=0.42]{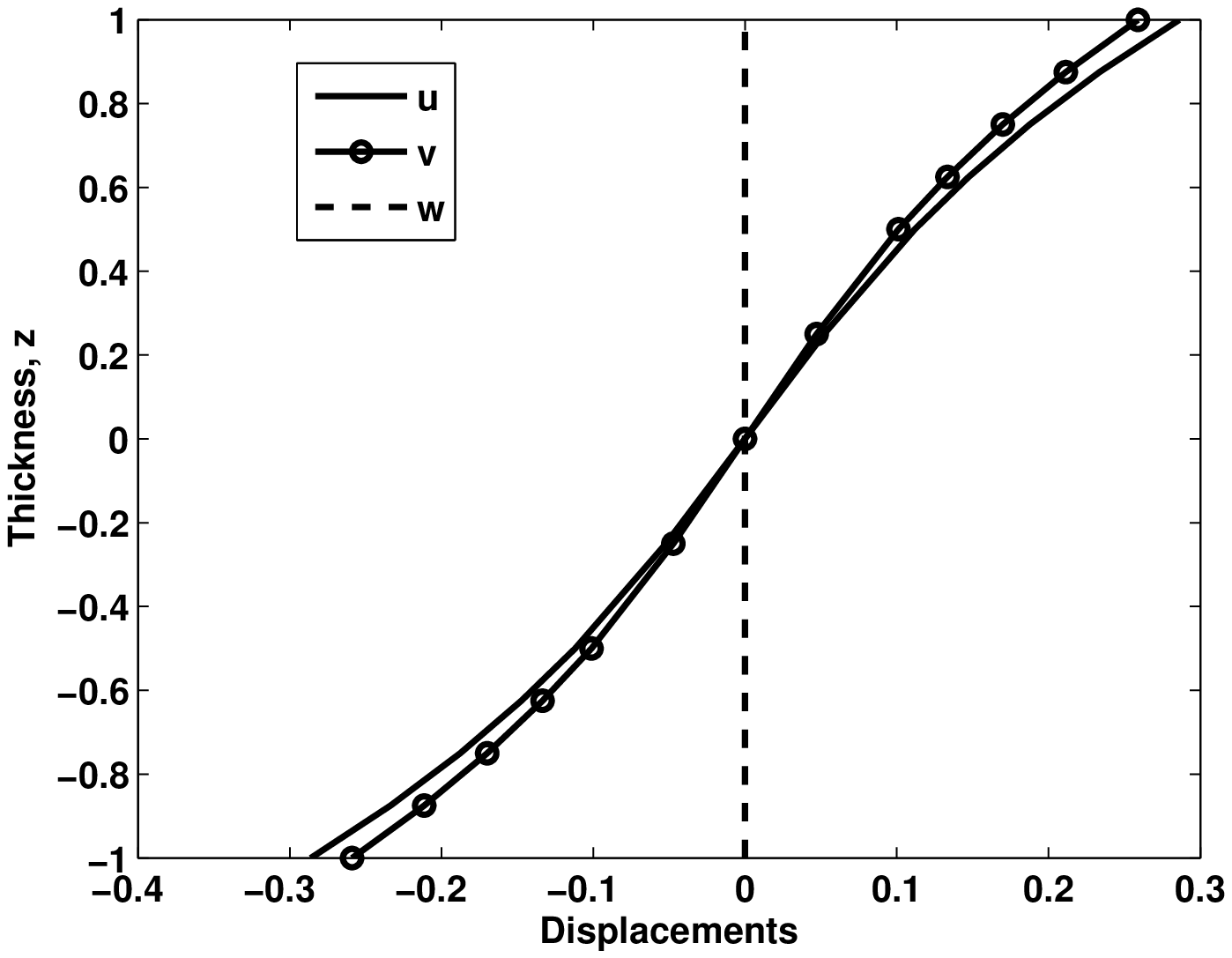}}
\caption{Deflected shapes $u_i(a/4,a/4,z),~i=x,y,z$ of the six modes for the square plates with simply supported edges for 1-2-1 Type A FGM plate with gradient index, $n =$ 1 and thickness, $a/h=$ 5.}
\label{fig:fgm121_quatA}
\end{figure}

\frefs{fig:fgm121_centerA} - \ref{fig:fgm121_quatA} shows the relative displacements through the thickness of the six modes using HSDT13. The plate is simply supported square 1-2-1 FGM plate with gradient index $n=$ 2 and $a/h=$ 5. The displacements $(u,v,w)$ are plotted along the lines $(a/2,b/2,z)$ and $(a/4,b/4,z)$, where $-h/2 \le z \le h/2$. It can be seen that, in flexural mode 1, the transverse displacement $w$ is not uniform at the chosen locations, thus, exhibiting the existence of normal stresses in the thickness direction. In flexural modes 2 - 6, the deflected shape retains the thickness at $(a/2,b/2,z)$, whereas the thickness of the plate is compressed at the other location, i.e., at $(a/4,b/4,z)$. It can be seen that the variation of the in-plane displacements $(u,v)$, in general, linear or non-linear in some higher-order models, irrespective of the types of modes.

\section{Conclusion}
\label{conclu}
FGM sandwich plate bending and free vibration analyses are carried out considering various parameters such as the material gradient index, the sandwich type and the thickness ratio. Different plate models are employed in predicting the physical behaviour and their through thickness variations in the plate. The accuracy and the effectiveness of the higher-order models (HSDT13 and HSDT11) over the lower-order theories have been demonstrated considering problems for which analytical/numerical results are available in the literature. It may be concluded that the HSDT13 model predict accurate results for any type of sandwich construction whereas other type of models depend on the type of sandwich plate and loading situations.

\bibliographystyle{elsarticle-num}
\bibliography{highorder}

\end{document}